\documentclass{au}

\newcommand{\<}{\kern.0833em}
\renewcommand{\[}{\begin{equation}\textstyle}
\renewcommand{\]}{\end{equation}}
\def\r#1{{\rm #1}}

\newtheorem{theorem}{Theorem}
\newtheorem{lemma}[theorem]{Lemma}
\newtheorem{corollary}[theorem]{Corollary}
\newtheorem{proposition}[theorem]{Proposition}
\newtheorem{definition}[theorem]{Definition}
\newtheorem{exercise}[theorem]{Exercise}

\newcommand{\Conv}{\mathrm{Conv}}
\newcommand{\RelConv}{\mathrm{RelConv}}
\newcommand{\R}{\mathbb R}
\newcommand{\V}{\mathbf V}
\newcommand{\odual}{{^\diamond}}
\newcommand{\Lat}{\mathbf{Lattice}}
\newcommand{\Latm}{\Lat^{m\text{-}\mathrm{pt}}}
\newcommand{\ch}{\mathrm{c.h.}}
\newcommand{\op}{^{\mathrm{op}}}
\newcommand{\inr}{\mathrm{int}}
\newcommand{\cl}{\mathrm{cl}}
\newcommand{\cp}{\mathrm{cpct}}
\newcommand{\ob}{\mathrm{o.bdd.}}

\begin{document}

\title{On lattices of convex sets in $\R^n$}
\thanks{
\protect
\protect\\
\indent
This preprint is readable online at
http://math.berkeley.edu/%
{$\!\sim$}gbergman/papers/\protect\linebreak[0]cnvx\<\_\,lat.%
\protect\linebreak[0]\{tex,\protect\linebreak[0]ps,dvi\},
and at arXiv:math.MG/0409288\,.
The versions on my website may be updated more often
than the arXiv copy.}

\author{George M. Bergman}
\address{University of California, Berkeley, CA 94720-3840, USA}
\email{gbergman@math.berkeley.edu}

\dedicatory{To the memory of Ivan Rival}

\begin{abstract}
Properties of several sorts of lattices of convex subsets
of $\R^n$ are examined.
The lattice of convex sets containing the
origin turns out, for $n>1,$ to satisfy a set of identities
strictly between those of the lattice of all convex subsets
of $\R^n$ and the lattice of all convex subsets of $\R^{n-1}.$
The lattices of arbitrary, of open bounded, and of compact convex
sets in $\R^n$ all satisfy the same identities, but
the last of these is join-semidistributive, while
for $n>1$ the first two are not.
The lattice of relatively convex subsets of a fixed
set $S\subseteq\R^n$ satisfies some, but in general not all of the
identities of the lattice of ``genuine'' convex subsets of $\R^n.$
\end{abstract}

\maketitle

\section{Notation, conventions, remarks.}\label{NHintro}

For $S$ a subset of $\R^n,$ the convex hull of $S$ will be denoted
\[\ch(S)\ =\ \{\Sigma_{i=1}^m\ \lambda_i\,p_i\mid m \geq 1,\ p_i\in S,
\ \lambda_i\in[0,1],\ \Sigma\,\lambda_i=1\}.\label{x.02}\]
When $S$ is written as a list of elements ``$\!\{\,...\,\},\!$'' we will
generally simplify ``$\!\ch(\{\,...\,\})\!$'' to
``$\!\ch(\,...\,)\!$''.

$\Conv(\R^n)$ will denote the lattice of
all convex subsets of $\R^n;$ its lattice operations are
\[\label{x.04} x\wedge y\ =\ x\cap y,\qquad x\vee y\ =\ \ch(x\cup y).\]
In any lattice, if a finite family of elements $y_i$ has
been specified, then an expression such
as $\bigwedge_i\,y_i$ will denote the
meet over the full range of the index $i,$ and similarly for joins.
Likewise, if we write something
like $\bigwedge_{j\neq i}\ y_j$ where $i$ has
been quantified outside this expression, then the meet will be over
all values of $j$ in the indexing family other than $i.$
If $L$ is any lattice and $x$ an element of $L,$ or, more generally,
of an overlattice of $L,$ we define the sublattice
\[\label{x.11} L_{\geq x}\ =\ \{y\in L\mid y\geq x\}.\]
In particular, $\Conv(\R^n)_{\geq\{0\}}$ is the
lattice of those convex subsets of $\R^n$ that contain the origin.

If $A$ is a class of lattices, $\V(A)$ will
denote the variety of lattices generated by $A,$ that is,
the class of lattices satisfying all identities (in the binary
operations $\wedge$ and $\vee)$ that hold in all lattices in $A.$
Again, for $A$ given as a list, we will
abbreviate $\V(\{\,...\,\})$ to $\V(\,...\,).$
(Most often $A$ will be a singleton $\{L\},$ so that
we will write $\V(L).)$

Given sets $x$ and $y,$ we will write $x-y$ for their set-theoretic
difference, $\{p\mid p\in x,\ p\notin y\}.$

Though I am not an expert either in lattice theory or in convex sets,
I know more about the former subject than the latter; hence,
I may more often state explicitly facts known to every worker in
convex sets than those known to every lattice-theorist.
I hope this note will nevertheless be of interest to people in
both fields.
I have no present plans of carrying these investigations
further; others are welcome to do so.

Since obtaining the main results of this paper, I have learned that
many of them were already in the literature, and have added references;
thus, this is now a hybrid research/survey paper.
I am grateful to
Kira Adaricheva, J.\,B.\,Nation,
Marina Semenova, Fred Wehrung,
and the referee, for corrections, information on the literature,
and many other helpful comments.

Whereas this note looks at conditions satisfied {\em universally}
in various lattices of convex sets, the papers \cite{KA.CvG},
\cite{FW+MS}, and others cited there study sufficient conditions
for lattices to be embeddable in such lattices, in other words,
{\em existential} properties of such lattices.
(This note includes one result of that type, in \S\ref{NHrlmsc}.)

\section{$n$\!-Distributivity.}\label{NHDn}

The varieties of lattices we will be examining in the first few
sections are those listed in
\begin{lemma}\label{X.13}
Each lattice in the sequence $\Conv(\R^0)_{\geq\{0\}},$ $\Conv(\R^0),$
$\Conv(\R^1)_{\geq\{0\}},$ $\dots,$ $\Conv(\R^n)_{\geq\{0\}},$
$\Conv(\R^n),\ \dots$ is embeddable in the next.
Hence
\[\label{x.14}
\V(\Conv(\R^0)_{\geq\{0\}})\subseteq\V(\Conv(\R^0))
\subseteq\V(\Conv(\R^1)_{\geq\{0\}})
\subseteq\V(\Conv(\R^1))\]
$$\textstyle \hspace{6em}\subseteq\dots
\subseteq \V(\Conv(\R^n)_{\geq\{0\}})
\subseteq \V(\Conv(\R^n))\subseteq\dots\ .\nonumber$$
\end{lemma}

\begin{proof}
On the one hand, $\Conv(\R^n)_{\geq\{0\}}$ is a sublattice
of $\Conv(\R^n);$ on the other, one can embed $\Conv(\R^n)$ in
$\Conv(\R^{n+1})_{\geq\{0\}}$ by sending each $x\in\Conv(\R^n)$ to the
cone $\ch(\{0\}\cup \{(p_1,\dots,p_n,1) \mid (p_1,\dots,p_n)\in x\}).$
The inclusions~(\ref{x.14}) follow from these embeddings.\end{proof}

Clearly $\Conv(\R^0)_{\geq\{0\}}$ is a trivial
lattice with unique element $\{0\};$ on the other hand, $\Conv(\R^0)$ is
a two-element lattice, so $\V(\Conv(\R^0))$ is
the variety of distributive lattices.
Hence the first inclusion of~(\ref{x.14}) is strict.
We shall see below that the next inclusion is an equality, while
all subsequent inclusions are again strict.

Let us set up notation for some lattice identities.
\begin{definition}[after Huhn \cite{AH.n-d}]\label{X.20}
For each positive integer $n,$ we shall denote
by $D_n(x,y_1,\dots,y_{n+1})$ \r(for ``$n$\!-distributivity''\/\r) the
lattice-relation in $n{+}2$ arguments $x,\,y_1,\,\dots,\,y_{n+1},$
\[\label{x.21}
x\wedge(\bigvee_i\,y_i)\ =\ \bigvee_i\ %
(x\wedge\bigvee_{j\neq i}\,y_j).\]
We shall say that a lattice $L$ satisfies ``the
identity $D_n\!\!$'' if\/~{\rm(\ref{x.21})} holds for
all $x,\,y_1,\linebreak[0]\,\dots,\,y_{n+1}\in L.$
\end{definition}

Thus, $D_1$ is the ordinary distributivity identity.

The use of the subscript $n$ for an identity
in $n{+}2$ variables which is symmetric
in $n{+}1$ of these may seem confusing;
a useful mnemonic is that $D_n$ is the identity
that allows one to ``reduce meets of $x$ with larger joins to
meets of $x$ with \!$n$\!-fold joins''.
An additional occasion for confusion will arise when we see
that the first of these identities to be satisfied
by $\Conv(\R^n)$ is not $D_n,$ but $D_{n+1}.$
This will be a consequence of the fact that an \!$n$\!-dimensional
simplex has $n{+}1$ vertices.

Note that the left-hand side of~(\ref{x.21}) is $\geq$ the
right-hand side for any family of elements of any lattice, since each of
the $n{+}1$ terms in the outer join on the right is
majorized by the left-hand side.
So to verify any instance of (\ref{x.21}) it suffices to
prove ``$\leq\!$''.

The pioneering work on identities satisfied by lattices
of convex sets was done by A.\,P.\,Huhn \cite{AH.Dop}, \cite{AH.n-d}.
The results in this and the next three sections will extend Huhn's
by approximately doubling both the family of lattices and the family of
identities considered, and formalizing some general techniques.
Huhn's results will be recovered along with our new ones.

The key to Huhn's and our results on $D_n$ (and some related identities)
is the following standard result in the theory of convex sets.
Strictly speaking, it is the first sentence below that is
Carath\'eodory's theorem, while the second is a well-known refinement
thereof \cite[p.431, line 4]{JE.HRC}.
Intuitively, that second sentence says that from an arbitrary
point $p_0\in S\subseteq \R^n,$ we can
``see'' any other point of $\ch(S)$ against a background of
(or embedded in) some \mbox{\!$(n{-}1)$\!-simplex} with
vertices in $S.$\bigskip
\\
{\bf Carath\'eodory's Theorem.} {\it
If $S$ is a subset of $\R^n,$ then each element $q\in\ch(S)$ belongs
to $\ch(p_0,\dots,p_n)$ for some $n{+}1$ points $p_0,\dots,p_n\in S.$
Moreover,
$p_0$ can be taken to be any pre-specified element of $S.$}\qed\medskip

In each paragraph of the next lemma, it is the first
assertion that is due to Huhn.

\begin{lemma}[cf.\ Huhn \cite{AH.n-d}]\label{X.23}
For every natural number $n$ and \!$(n{+}3)$\!-tuple
of convex sets $x, y_1, \dots, y_{n+2}\in\Conv(\R^n),$ the
relation $D_{n+1}(x,y_1,\dots,y_{n+2})$ holds.
Moreover, if $n\,{>}\,0$ and \r(at least\/\r) some two of
the sets $y_1,\dots,y_{n+1}$ have
nonempty intersection, then $D_n(x,y_1,\dots,y_{n+1})$ holds.

Hence $\Conv(\R^n)$ satisfies the identity $D_{n+1},$ and
for $n>0,$ $\Conv(\R^n)_{\geq\{0\}}$ satisfies the identity~$D_n.$
\end{lemma}
\begin{proof}
We shall prove the assertions of the first paragraph, which clearly
imply those of the second.

To get the first assertion,
let $x,\,y_1,\dots,y_{n+2}\in\Conv(\R^n),$ and let $p$ be a
point of the convex set described by the left-hand side of~(\ref{x.21}).
Then $p$ belongs to both $x$ and $\bigvee_i\,y_i=\ch(\bigcup_i\, y_i).$
Carath\'eodory's Theorem now says that $p$ belongs to the convex hull
of $n{+}1$ points of $\bigcup_i\,y_i,$ hence to the join of at most
$n{+}1$ of the $y_i,$ so it belongs to one of the terms on the
right-hand side of~(\ref{x.21}), hence to their join, as required.

To prove the second assertion,
let $x,\ y_1,\,\dots,\,y_{n+1}\in\Conv(\R^n),$ and suppose that two
of $y_1,\dots,y_{n+1}$ have a point $p_0$ in common.
Given $p$ in the left-hand side of the desired instance
of~(\ref{x.21}), we see from the last sentence of our statement
of Carath\'eodory's Theorem that $p$ will belong to the convex
hull of $p_0$ and some $n$ other points of $\bigcup_i\,y_i.$
The latter points will lie in the union of some $n$ of
the $y$'s, and since those $n$ of
our $n{+}1$ sets leave out only one, they do
not leave out both of the sets known to contain $p_0.$
So that family of $n$ $y$'s contains the $n{+}1$ points whose convex
hull is known to contain $p,$ and the conclusion
follows as before.\end{proof}

Note that for all $p\in\R^n,$ we
have $\Conv(\R^n)_{\geq\{p\}}\cong\Conv(\R^n)_{\geq\{0\}};$ moreover,
if $p\in S\subseteq\R^n,$ then
$\Conv(\R^n)_{\geq S}\subseteq \Conv(\R^n)_{\geq\{p\}}.$
Thus, any identities proved for the
lattice $\Conv(\R^n)_{\geq\{0\}}$ will hold
in $\Conv(\R^n)_{\geq S}$ for every nonempty set~$S.$

One may ask whether, in the first assertion of the above lemma, we
have failed to use the full strength of Carath\'eodory's Theorem.
That theorem says that the convex hull of any family
of $N>n{+}1$ points is the union of the hulls of its
\!$(n{+}1)$\!-element subfamilies, hence for
each such $N,$ we may deduce an identity like $D_{n+1},$ but
with the join of $N$ rather than $n{+}2$ convex sets $y_i$ on the left,
and expressions involving all the \!$(n{+}1)$\!-fold subjoins
thereof on the right.
Our $D_{n+1}$ is the case $N=n{+}2.$

But in fact, the identities so obtained are all equivalent to $D_{n+1}.$
To derive them from it, note first that if in $D_{n+1}$ we
substitute for $y_{n+2}$ the expression $y_{n+2}\vee y_{n+3},$ then the
left-hand side becomes $x\wedge(\bigvee_{i=1}^{n+3}y_i),$ while
on the right, some of the joins involve $n{+}1$ of
the $y$'s and others involve $n{+}2.$
If we again apply $D_{n+1}$ to the latter joins, we get
precisely the $N=n{+}3$ case of the class of identities discussed above.
The identities with still larger $N$ are gotten by repeating
this argument.
Conversely, one can get $D_{n+1}$ from any of these identities by
substituting for the $y_i$ with $i>n{+}2$ repetitions
of $y_{n+2},$ and discarding from the outer join
on the right-hand side joinands majorized by others.

Let us also note that $D_{n+1}{\implies}D_{n+2}.$
Indeed, the identity with $n{+}3$ $y$'s
that we just showed equivalent to $D_{n+1}$ has
the same left-hand side as $D_{n+2},$ while the right-hand
side of $D_{n+2}$ can be seen to lie, in an arbitrary
lattice, between the two sides of that identity.
In particular, the identities obtained in
Lemma~\ref{X.23} are successively
weaker for larger $n,$ as is reasonable in view of~(\ref{x.14}).

The argument proving Lemma~\ref{X.23} (for simplicity let us limit
ourselves to the first assertion thereof) can be formulated in a
more general context, and the above observations on identities
allow us to obtain a converse in that context, which we record below,
though we shall not use it.
Recall that a closure operator on a set $X$ is called
``finitary'' (or ``algebraic'') if the closure of every
subset $S\subseteq X$ is
the union of the closures of the finite subsets of $S.$

\begin{lemma}[cf.\ \cite{LL.n-d}]\label{X.24}
Let $\cl$ be a finitary closure operator on a set $X,$ such that
every singleton subset of $X$ is closed, or more generally, such that
the closure of every singleton is join-irreducible;
and let $n$ be a positive integer.
Then the lattice of closed subsets
of $X$ under $\cl$ satisfies $D_n$ if and only
if $\cl$ has the \mbox{``$n$\!-Carath\'eodory}
property'' that the closure of every set is the union of the
closures of its ${\leq}\<n\<\!$-element subsets.
\end{lemma}
\begin{proof}
``If'' is shown exactly as in the proof of Lemma~\ref{X.23}.

Conversely, suppose the lattice of \!$\cl$\!-closed subsets
of $X$ satisfies $D_n,$ and
let $p\in\cl(S)$ for some $S\subseteq X.$
By the assumption that $\cl$ is
finitary, $p\in\cl(q_1,\dots,q_N)$ for some $q_1,\dots,q_N\in S.$
If $N\leq n,$ we are done; if not, let us
rewrite the condition $p\in\cl(q_1,\dots,q_N)$ as a
lattice relation, $\cl(p)=\cl(p)\wedge (\bigvee_i\cl(q_i)).$
By the preceding discussion, $D_n$ implies the identity
in $N$ variables which, when applied to the right-hand side of
the above relation, turns the relation into
$$\textstyle\cl(p)\ =\ \bigvee_{|I|=n}\ (\cl(p)\,\wedge\,
(\bigvee_{i\in I}\ \cl(q_i))),$$
where the subscript to the outer join means that $I$ ranges
over all \!$n$\!-element subsets of $\{1,\dots,N\}.$
Now by assumption $\cl(p)$ is join-irreducible, hence
it equals one of the joinands on the right,
showing $p$ to be in the
closure of some set of $n$ $q$'s, as required.\end{proof}

\section{Tools for studying related identities.}\label{NHudud}

We also appear to have used less than the
full force of the middle sentence of Lemma~\ref{X.23} (the stronger
relation holding when at least two of the $y$'s have nonempty
intersection) in getting the second assertion of the last
sentence of that lemma (the stronger identity
for $\Conv(\R^n)_{\geq\{0\}}),$ since
the latter assertion concerns the
case where not just two, but all the $y_i$ (and
also $x)$ have a point in common.
There is no evident way to take advantage of the weaker hypothesis
of said middle sentence when working in the
lattice $\Conv(\R^n)_{\geq\{0\}};$ but might we
be able to use it to prove some new identity
for $\Conv(\R^n);$ for example, by finding a
relation that holds, on the one hand, whenever a family of elements
$x,y_1,\dots,y_{n+1}$ satisfies $D_n,$ and also, for some trivial
reason, whenever $y_1\wedge y_2=\emptyset$?

One relation with these properties can be obtained
by taking the meet of each side of $D_n$ with $y_1\wedge y_2.$
Unfortunately, this turns out to be the trivial identity: both sides
simplify to $x\wedge y_1\wedge y_2$ in any lattice.
However, we can circumvent this by first taking the join of
both sides of $D_n$ with a new indeterminate $z,$ and only then
taking meets with $y_1\wedge y_2.$

We shall in fact see that this provides what Lemma~\ref{X.23}
failed to: an identity holding in $\Conv(\R^n)$ but
not in $\Conv(\R^{n+1})_{\geq\{0\}}.$
However, the verification of an example showing the failure
of this identity in the larger lattice (and of a similar example
we will need later) is messy if done entirely ``by hand''; so we shall
establish in this section some general criteria for certain sorts of
inequalities to be strict.

The last assertions of the next two lemmas clearly imply,
for $\Conv(\R^n)$ and $\Conv(\R^{n+1})_{\geq\{0\}}$ respectively, that
if we are given expressions $u,$ $v$ and $w$ in some lattice
indeterminates such that the inequality $u\geq v$ is known
to hold identically in our lattice, then we can write down another
inequality holding identically among expressions in a
slightly larger set of indeterminates, for which equality will
hold whenever either $u$ equals $v,$ or $w$ equals
the least lattice-element (i.e., $\emptyset,$ respectively $\{0\}),$ but
which will fail in all other cases.
This is what is logically called for by the program sketched above.
However the earlier parts of these lemmas give some simpler
inequalities for which the same is true
if the hypotheses hold ``in a sufficiently strong way'', and
it will turn out that in the applications where we need to show failure
of an identity, we will be able to use these simpler formulas.

In the statements of these lemmas,
for $p, q\in\R^n,$ ``the ray drawn from $q$ through $p\!$'' will
mean $\{\lambda\<p + (1-\lambda)q\mid 0 \leq \lambda < \infty\},$ even
in the degenerate case $p=q,$ where this set is the singleton $\{p\}.$
\begin{lemma}\label{X.28}
Let $n$ be any natural number, and let $u \geq v$ and $w$ be three
elements of $\Conv(\R^n).$
Then\smallskip
\\
{\rm(i)} The following conditions are equivalent:

{\rm(a)} There exists $z\in\Conv(\R^n)$ such
that $(u\vee z)\wedge w > (v\vee z)\wedge w.$

{\rm(b)} There exist a point $p\in u$ and
a point $q\in w$ such that the ray drawn
from $q$ through $p$ contains no point of $v.$\smallskip
\\
{\rm(ii)} The following conditions are equivalent:

{\rm(a)} There exist $z, z'\in\Conv(\R^n)$ such
that $((u\wedge z')\vee z)\wedge w > ((v\wedge z')\vee z)\wedge w.$

{\rm(b)} $u>v,$ and $w$ is nonempty.
\end{lemma}
\begin{proof}
In view of the hypothesis $u \geq v,$ the
relation ``$\geq\!$'' always holds in the inequalities
of~(i)(a) and~(ii)(a), so in each case, strict inequality
is equivalent to the existence of a point $q$ belonging to the
left-hand side but not to the right-hand side.

Suppose, first, that (i)(a) holds for some $z;$ thus we get
a point $q$ belonging to $u\vee z$ and to $w$ but not to $v\vee z.$
Note that the latter condition implies that $q$ does not lie in $z.$
If $q$ lies in $u,$ then~(i)(b) is
satisfied with $p=q,$ so assume $q\notin u.$
Hence, being in the convex hull of $u$ and $z$ but
in neither set, $q$ must lie on the
line-segment $\ch(p,r)$ for some $p\in u,\ r\in z:$

\begin{center}
\begin{picture}(100,50)(0,-25)
\thinlines
\put(-2,0){\line(-1,0){10}}
\put(2,0){\line(1,0){56}}
\put(62,0){\line(1,0){56}}
\put(122,0){\line(1,0){10}}

\put(0,0){\circle{3}}
\put(60,0){\circle{3}}
\put(120,0){\circle{3}}

\put(0,8){\makebox(0,0)[c]{$p$}}
\put(0,-8){\makebox(0,0)[c]{$\in u$}}
\put(60,8){\makebox(0,0)[c]{$q$}}
\put(60,-8){\makebox(0,0)[c]{$\in(u\vee z)\wedge w$}}
\put(60,-18){\makebox(0,0)[c]{$\notin v\vee z$}}
\put(120,8){\makebox(0,0)[c]{$r$}}
\put(120,-8){\makebox(0,0)[c]{$\in z$}}

\end{picture}
\end{center}

\noindent
Now if a point $s\in v$ lay on the ray
drawn from $q$ through $p,$ we would
have $q\in\ch(s,r) \subseteq v\vee z,$ contradicting our assumption
that $q\notin v\vee z.$
This proves~(b).

Conversely, if we are given $p$ and $q$ as in (i)(b), take $z=\{2q-p\}:$

\begin{center}
\begin{picture}(100,50)(0,-25)
\thinlines
\put(-2,0){\line(-1,0){10}}
\put(2,0){\line(1,0){56}}
\put(62,0){\line(1,0){56}}
\put(122,0){\line(1,0){10}}

\put(0,0){\circle{3}}
\put(60,0){\circle{3}}
\put(120,0){\circle{3}}

\put(0,8){\makebox(0,0)[c]{$p$}}
\put(0,-8){\makebox(0,0)[c]{$\in u$}}
\put(60,8){\makebox(0,0)[c]{$q$}}
\put(60,-8){\makebox(0,0)[c]{$\in w$}}
\put(120,8){\makebox(0,0)[c]{$2q-p$}}
\put(120,-8){\makebox(0,0)[c]{$\in z$}}

\end{picture}
\end{center}

\noindent
Thus $q\in(u\vee z)\wedge w,$ but we claim that $q\notin v\vee z.$
Indeed, if $q$ belonged to this set,
then $v$ would have to meet the ray from the unique
point $2q-p$ of $z$ through $q$ on
the other side of $q;$ i.e., it would meet the ray drawn
from $q$ through $p,$ contradicting our choice of $p$ and $q$ as in~(b).
Hence $q$ belongs to the left-hand but not the right-hand
side of the inequality of~(a), as required.

Note that~(i)(b), and hence~(i)(a), holds
whenever $u$ and $w$ are nonempty and $v$ is empty.

Now in the situation of (ii)(b), if we take
any $p\in u-v$ and let $z'=\{p\},$ then $u\wedge z'$
and $v\wedge z'$ are respectively nonempty and empty,
so applying the preceding observation with these two sets in
the roles of $u$ and $v,$ we get~(ii)(a).
The reverse implication is trivial.\end{proof}

The result we shall prove
for $\Conv(\R^n)_{\geq\{0\}}$ is similar.
Indeed, in the lemma below, part~(i) is
exactly as in the preceding lemma; but~(ii)
becomes two statements,~(ii) and~(iii), the former having the same
``(a)'' as in~(ii) above but a stronger ``(b)'',
the latter a weaker ``(a)'' and the same ``(b)'' as above.
Note also the restriction on $n$ (only needed for~(iii)).
\begin{lemma}\label{X.29}
Let $n>1,$ and let $u \geq v$ and $w$ be three
elements of $\Conv(\R^n)_{\geq\{0\}}.$
Then\smallskip
\\
{\rm(i)} The following conditions are equivalent:

{\rm(a)} There exists $z\in\Conv(\R^n)_{\geq\{0\}}$ such
that $(u\vee z)\wedge w\ >\ (v\vee z)\wedge w.$

{\rm(b)} There exist a point $p\in u$ and
a point $q\in w$ such that the ray drawn
from $q$ through $p$ contains no point of~$v.$\smallskip
\\
{\rm(ii)} The following conditions are equivalent:

{\rm(a)} There exist $z,z'\in\Conv(\R^n)_{\geq\{0\}}$ such
that
$$\textstyle ((u\wedge z')\vee z)\wedge
w\ >\ ((v\wedge z')\vee z)\wedge w.$$

{\rm(b)} There exist a point $p\in u$ and a
point $q\in w$ such that the ray drawn from $q$ through $p$ contains
no point of $v\wedge \ch(0,p).$\smallskip
\\
{\rm(iii)} The following conditions are equivalent:

{\rm(a)} There exist $z, z', z'', z'''\in\Conv(\R^n)_{\geq\{0\}}$ such
that
$$\textstyle((((u\wedge z''')\vee z'')\wedge z')\vee z)\wedge
w\ >\ ((((v\wedge z''')\vee z'')\wedge z')\vee z)\wedge w.$$

{\rm(b)} $u > v,$ and $w \neq \{0\}.$

\end{lemma}
\begin{proof}
(i)(a)${\implies}$(i)(b) holds by the preceding lemma.
In proving the reverse implication, we cannot
set $z=\{2q-p\}$ as we did there, so let $z=\ch(0,2q-p).$
As before, we have $q\in(u\vee z)\wedge w$ and
need to show $q\notin v\vee z.$
If the contrary were true, then $q$ would be a convex linear
combination of $0,$ $2q-p,$ and a point $r\in v.$
This can be rewritten as a convex linear combination of
$2q-p$ with a convex linear combination of $0$ and $r;$ but the latter
combination would also be a point of $v,$ and, as
in the previous proof, would lie on the ray drawn
from $q$ through $p,$ contradicting our choice
of $p$ and $q$ as in~(b).

Turning to (ii), if~(ii)(a) holds then we can
apply~(i)(a)${\implies}$(i)(b)
with $u\wedge z'$ and $v\wedge z'$ in place
of $u$ and $v,$ and the
resulting $p$ and $q$ will satisfy (ii)(b) (since
$v\wedge \ch(0,p) \subseteq v\wedge z').$
Inversely, if (ii)(b) holds, take $z'=\ch(0,p),$ and
apply~(i)(b)${\implies}$(i)(a)
with $u\wedge z'$ and $v\wedge z'$ in place of $u$ and $v.$

In statement (iii), it is clear that (a) implies (b).
To prove the converse, let us assume~(iii)(b), and consider two
cases, according to whether the stronger statement (ii)(b) holds.
If it does, we get the inequality of (ii)(a), from which
we can immediately get that of (iii)(a)
by choosing $z'''$ and $z''$ to ``have no effect'' (e.g., by taking them
to be $u$ and $v$ respectively).
On the other hand, if (ii)(b) fails while (iii)(b) holds, it is easy
to see that all elements of $w$ and $u-v$ must
lie on a common line $x$ through~$0.$
In that case, we want to
use $z'''$ and $z''$ to
``perturb'' $u$ and $v,$ so that the modified $u$ has
points off the line $x,$ while being careful to preserve the
property that $u$ is strictly larger than $v.$
To do this, we begin by taking any
point $p\in u-v,$ letting $z'''=\ch(0,p),$ and
noting that the intersections of $u$ and $v$ with
this segment are still distinct.
Now taking any point $r$ not on the line $x$ (it is for
this that we need $n>1),$ and letting $z''=\ch(0,r),$ we see
that $(u\wedge z''')\vee z''$ and $(v\wedge z''')\vee z''$ remain
distinct, and that their difference now has points off $x.$
Hence (ii)(b) holds with these sets in the roles
of $u$ and $v,$ and the implication (ii)(b)${\implies}$(ii)(a) gives
the $z$ and $z'$ needed for~(iii)(a).\end{proof}

One can get criteria similar to those of the preceding lemmas
for other conditions.
At the trivial end, given $u\in\Conv(\R^n),$ a condition
for $u$ to be nonempty is that there exist a $z$ such
that $u\vee z>z,$ and likewise the condition for at least
one of two elements $u$ and $v$ to be
nonempty is that their join have this property.
The exercise below offers, for the diversion of the
interested reader, some less trivial cases.
(We will not use the results of this exercise.)

\begin{exercise}\label{X.31}
{\rm(i)} Find an inequality in elements $u,\,v$ and
one or more additional lattice
variables $z, \dots,$ which holds identically in lattices,
and which, for any $u, v\in\Conv(\R^n),$ is
strict for some values of the additional variables if and
only if {\em both} $u$ and $v$ are nonempty.\smallskip
\\
{\rm(ii)} Find an inequality in elements $u,\,v,\,w,\,x$ and
additional variables which holds in any lattice
when $u \geq v,$ $w \geq x,$ and which, for any
such $u, v, w, x\in\Conv(\R^n),$ is strict for some values of the
additional variables if and only if $u>v$ and $w>x.$\smallskip
\\
{\rm(iii)} Same as {\rm(ii)}, but with ``$u>v$ and $w>x\!$''
replaced by ``$u>v$ or $w>x\!$''.\smallskip
\\
{\rm(iv)-(vi)} Like {\rm(i)-(iii)}, but for
$\Conv(\R^n)_{\geq\{0\}},$ and with ``nonzero'' in place of ``nonempty''
in~{\rm(i)}.\smallskip
\\
{\rm(vii)} In Lemma~\ref{X.29}, conditions
{\rm(i)(a)}, {\rm(ii)(a)} and {\rm(iii)(a)} involved 1, 2 and 4
$z\!$'s respectively; thus the condition that there exist three
elements $z,z',z''\in\Conv(\R^n)_{\geq\{0\}}$ such that
$$\textstyle (((u\vee z'')\wedge z')\vee z)\wedge w>
(((v\vee z'')\wedge z')\vee z)\wedge w$$
was skipped.
Show by example that this condition is not equivalent
to~{\rm(iii)(b)}.\smallskip
\\
{\rm(viii)} Suppose $u \geq v$ and $w$ are
elements of $\Conv(\R^n),$ and consider the conditions dual
to~{\rm(i)(a)-(iii)(a)} of Lemma~\ref{X.29}, and to the ``skipped''
condition:\smallskip

{\rm(a)} There exists $z\in\Conv(\R^n)$ such that
$(u\wedge z)\vee w>(v\wedge z)\vee w.$

{\rm(a$\!{}'$)} There exist $z, z'\in\Conv(\R^n)$ such
that $((u\vee z')\wedge z)\vee w > ((v\vee z')\wedge z)\vee w.$

{\rm(a$\!{}''$)} There exist $z,z',z''\in\Conv(\R^n)$ such that
$$\textstyle(((u\wedge z'')\vee z')\wedge z)\vee
w\ >\ (((v\wedge z'') \vee z')\wedge z)\vee w.$$

{\rm(a$\!{}'''$)} There exist $z,z',z'',z'''\in\Conv(\R^n)$ such that
$$\textstyle((((u\vee z''')\wedge z'') \vee z')\wedge z)\vee
w\ >\ ((((v\vee z''')\wedge z'')\vee z')\wedge z)\vee w.$$

\noindent
Which of these, if any, are equivalent, for all
such $u,$ $v$ and $w,$ to the condition\smallskip

{\rm(b)} $u>v,$ and $w\neq \R^n$?
\end{exercise}

Now, back to business.

\section{Identities distinguishing our chain of lattices.}

Given $n>0,$ let us write $((D_n)\vee z)\wedge y_1\wedge y_2$ for
the equation in $n{+}3$ variables $x, y_1, \dots, y_{n+1}, z$ obtained
by applying the operator $((-)\vee z)\wedge y_1\wedge y_2$ to both
sides of the relation $D_n.$
We can now prove
\begin{theorem}\label{X.25}
For each positive integer $n,$ $\Conv(\R^n)$
satisfies $((D_n)\vee z)\wedge y_1\wedge y_2$ but
not $D_n,$ while $\Conv(\R^n)_{\geq\{0\}}$ satisfies
$D_n$ but, if $n>1,$ not $((D_{n-1})\vee z)\wedge y_1\wedge y_2.$

Hence every inclusion in {\rm(\ref{x.14})} is strict except the second.
Equality holds at that step.
\end{theorem}
\begin{proof}
To see that equality holds at the second
inclusion of~(\ref{x.14}), note that both $\Conv(\R^0)$ and
$\Conv(\R^1)_{\geq\{0\}}$ are nontrivial lattices,
which by Lemma~\ref{X.23} satisfy $D_1,$ the distributive
identity, and that the variety of distributive lattices is known
to have no proper nontrivial subvarieties.

The positive assertions of the first paragraph of the theorem follow
from Lemma~\ref{X.23}, combined, in the case of
the first of these results,
with the observations of the first two paragraphs of \S\ref{NHudud}.
(Those observations are equivalent to the contrapositive of the
easy implications
(i)(a)${\implies}$(ii)(a)${\implies}$(ii)(b) of Lemma~\ref{X.28}.)
It remains to give examples showing the negative assertions.

To see that $\Conv(\R^n)$ does not
satisfy $D_n,$ let $q_1, \dots, q_{n+1}$ be the vertices of an
\!$n$\!-simplex in $\R^n$ and $p$ an
interior point of this simplex, and take
for the $y_i$ and $x$ the
singletons $\{q_i\}$ and $\{p\}$ respectively.
Then we see that the left-hand side
of $D_n(x,y_1,\dots,y_{n+1})$ gives $x,$ while all the joinands
on the right are empty, hence so is the right-hand side itself.

To show that $\Conv(\R^n)_{\geq\{0\}}$ does not
satisfy $((D_{n-1})\vee z)\wedge y_1\wedge y_2$ when $n>1,$ let
$P\subseteq \R^n$ be a hyperplane not passing through $0.$
The idea will be to mimic the preceding example within $P,$ then
replace the resulting singleton sets with the line-segments connecting
them with $0,$ slightly
enlarge $y_1,$ so that it has nonzero intersection
with $y_2,$ and finally apply part~(i) of Lemma~\ref{X.29}.

So let $q_1,\dots,q_n$ be the
vertices of an \!$(n{-}1)$\!-simplex in $P,$ and $p$ a
point in the relative interior of that simplex, and
let $x$ and the $y_i$ be the line
segments $\ch(0,p)$ and $\ch(0,q_i)$ respectively,
except for $y_1,$ which we
take to be $\ch(0,\<q_1,\<q_2/2).$
Note that all of these convex sets lie in the closed
half-space $H$ bounded by $P$ and containing~$0;$ hence
the intersection with $P$ of any lattice expression in
these convex sets can be computed as the corresponding lattice
expression in their intersections with $P.$
We see that these intersections are a configuration of
the form given in the preceding example, except that the
dimension is lower by~$1$
(note that the point ``$q_2/2\!$'' in
our definition of $y_1,$ not lying in $P,$ does
not affect the intersection $y_1\cap P).$
Hence when we evaluate the two sides of $D_{n-1}$ at these
elements, the left-hand side intersects $P$ in
the point $p,$ and is thus the whole
line-segment $x,$ while the right-hand side does not
meet $P,$ hence is a proper subsegment of $x.$

We can now deduce the failure of
$((D_{n-1})\vee z)\wedge y_1\wedge y_2$ from the
implication (i)(b)${\implies}$(i)(a) of Lemma~\ref{X.29},
using for $u$ and $v$ respectively the left and right sides
of $D_{n-1},$ and for $w$ the set $y_1\wedge y_2.$
To see that~(i)(b) holds for these sets, note
that $w$ is the line-segment $\ch(0,q_2/2).$
Since this lies in a different line through the origin
from $p,$ the ray drawn from any nonzero point $q$ of
that segment through $p$ does not
meet the line-segment $x$ in any point other than $p,$ so
in particular, it does not contain any point of $v,$ the
right-hand side of $D_{n-1},$ which we saw was a proper
subset of $x.$
\end{proof}

We remark that by alternately inserting joins and meets with more and
more variables $z^{(m)}$ into
``$\!((D_n)\,\dots\,) \wedge(y_1\wedge y_2)\!$''
one can get identities that, formally, are successively stronger
(though still all implied by $D_n),$ so that the statements that a
lattice does not satisfy these identities become successively weaker.
Thus, as a stronger version of the above theorem, we could have
stated that $\Conv(\R^n)$ satisfies such identities
with arbitrarily long strings of inserted terms,
while $\Conv(\R^n)_{\geq\{0\}}$ fails to
satisfy the particular one given above.
But for simplicity, I used just one identity to
distinguish the properties of these lattices.

The argument at the beginning of the preceding
section showing that $\Conv(\R^n)$ satisfies
$((D_n)\vee z)\wedge y_1\wedge y_2$ also
clearly shows that it satisfies the formally stronger identity
\[\label{x.27}
((D_n)\vee z)\,\wedge\,
(\bigvee_{i,j;\ i\,\neq\,j}\ y_i \wedge y_j).\]
With a little additional work one can get the still stronger identity:
\[\label{x.26}
((D_n)\vee z)\,\wedge\,\bigwedge_i\,(\bigvee_{j\neq i}\,y_j).\]
(Idea: If $\bigwedge_i\,(\bigvee_{j\neq i}\,y_j)$ is nonempty,
use a point thereof as the $p_0$ in the
second sentence of Carath\'eodory's Theorem.)

We remark that the fact that $\Conv(\R^n)$ satisfies
the identity $D_{n+1}$ says that
for any $n{+}2$ convex sets $y_i,$ the
union of the \!$(n{+}1)$\!-fold
joins $\bigvee_{j\neq i}\,y_j$ $(i=1,\dots,n{+}2)$ is
itself convex, while~(\ref{x.26})
says essentially that the same holds for the union of
the \!$n$\!-fold joins of $n{+}1$ convex sets,
if those unions have at least one point in common.

\section{Dual \!$n$\!-distributivity.}\label{NHDnop}

Huhn showed not only that $\Conv(\R^n)$ satisfies
$D_{n+1},$ but also that it satisfies the dual of that identity.
Let us write the dual of the relation $D_n(x,y_1,\dots,y_{n+1})$ as
\[\label{x.30}
D_n\op(x,y_1,\dots,y_{n+1}):\quad x \vee (\bigwedge_i\ y_i)\,=\,
\bigwedge_i\ (x\vee\bigwedge_{j\neq i}\,y_j).\]
Just as in $D_n$ the direction $\geq$ is automatic,
so $\leq$ is automatic in $D_n\op.$

Below, we will strengthen Huhn's result
that $D\op_{n+1}$ holds identically
in $\Conv(\R^n)$ by showing that it holds in the larger
lattice $\Conv(\R^{n+1})_{\geq\{0\}},$ and will
again use the method of \S\ref{NHudud} to manufacture
a related identity which holds in $\Conv(\R^n)$ but not
in $\Conv(\R^{n+1})_{\geq\{0\}}.$
Like Huhn, we start with\bigskip
\\
{\bf Helly's Theorem} (\cite{JE.HRC}, p.391){\bf .} {\it
Let $n \geq 0.$
If a finite family of convex subsets of $\R^n$ has the
property that every $n{+}1$ of them have nonempty intersection,
then the whole family has nonempty intersection.}\qed\medskip

We will also use the following observation.
\begin{lemma}\label{X.33}
Let $V$ be any real vector space.
Given a convex set $x$ and a point $p$ in $V,$ there
exists a convex set $w$ in $V$ such that for every
nonempty convex set $y,$ one has $p\in x\vee y$ if
and only if $y$ has nonempty intersection with $w.$
\end{lemma}
\begin{proof}
It is straightforward to check that a set $w$ with the required
property
(in fact, the unique such set) is the cone consisting of the union
of all rays from $p$ which meet the central reflection
of $x$ through $p.$
(If $p\in x$ then $w=V;$ in
visualizing the contrary case, it is convenient to assume
without loss of generality that $p=0.)$\end{proof}

We can now show that certain sorts of families of
convex sets satisfy the relation $D_n\op,$ from
which we will deduce our identities.
\begin{lemma}\label{X.35}
Let $n$ be a natural number, and $x,y_1,\dots,y_{n+1}$ be
$n{+}2$ elements of $\Conv(\R^n)$ such
that $y_1\wedge\dots\wedge y_{n+1}\neq \emptyset.$
Then $D_n\op(x,y_1,\dots,y_{n+1})$ holds.
\end{lemma}
\begin{proof}
Given
\[\label{x.36}
p\in\bigwedge_i\ (x\vee\bigwedge_{j\neq i}\,y_j)\]
we need to show that
\[\label{x.37}
p\in x \vee (\bigwedge_i\,y_i).\]
Let $w$ be the set determined by $x$ and $p$ as in Lemma~\ref{X.33}.
Since by~(\ref{x.36}), $p$ lies in each of the sets
$x\vee\bigwedge_{j\neq i} y_j,$ our choice of $w$ shows
that for each $i,$ $(\bigwedge _{j\neq i} y_j)\wedge w$ is nonempty.
These conditions together with the nonemptiness
of $y_1\wedge\dots\wedge y_{n+1}$ allow us to
apply Helly's Theorem to the $n{+}2$ convex
sets $y_1, \dots, y_{n+1}, w,$ and conclude
that $(\bigwedge_i y_i)\wedge w$ is nonempty, which
by choice of $w$ is equivalent to~(\ref{x.37}).\end{proof}

Now for each positive integer $n,$
let $(((D_n\op)\wedge z')\vee z) \wedge(\bigwedge_i y_i)$ denote
the relation in variables $x,$ $y_1,$ $\dots,$ $y_{n+1},$ $z,$ $z'$
obtained by applying the
operation $(((-)\wedge z')\vee z)\wedge (\bigwedge_i y_i)$ to
both sides of~(\ref{x.30}).
Then we have
\begin{theorem}[cf.\ Huhn \cite{AH.Dop}]\label{X.40}
For each positive integer $n,$ $\Conv(\R^n)$ satisfies the
identity $(((D_n\op)\wedge z') \vee z)\wedge(\bigwedge_i\,y_i)$ but
not $D_n\op,$ while $\Conv(\R^n)_{\geq\{0\}}$ satisfies the
identity $D_n\op$ but, if $n>1,$ not
$(((D\op_{n-1}) \wedge z')\vee z)\wedge (\bigwedge_i\,y_i).$
\end{theorem}
\begin{proof}
The positive assertions follow from Lemma~\ref{X.35}.
(Again, the reasoning that obtains the first identity from that lemma
can be considered, formally, an application of the contrapositive
of an easy direction in of one of our lemmas, in this case the
implication (ii)(a)${\implies}$(iii)(b) of Lemma~\ref{X.29}.)

To get a counterexample
to $D_n\op$ in $\Conv(\R^n),$ take for $x$ an \!$n$\!-simplex
in $\R^n,$ let $p$ be a point
outside $x,$ let $x'$ be the central reflection of $x$ through $p,$ and
let $y_1,\dots,y_{n+1}$ be the \!$(n{-}1)$\!-faces of $x'.$
Then $\bigwedge_i\,y_i$ is empty, so the
left-hand side of~(\ref{x.30}) is just $x,$ while each of the
intersections $\bigwedge_{j\neq i} y_j$ is a nonempty subset of $x',$ so
the right-hand side of~(\ref{x.30}) contains $p.$

As in the proof of Theorem~\ref{X.25}, we begin the counterexample
to our more elaborate identity
in $\Conv(\R^n)_{\geq\{0\}}$ by taking a
copy of our preceding example, for the next lower dimension, in
a hyperplane $P \subseteq \R^n$ not containing~$0.$
Let us write $p_0$ and $x_0,$ $x_0',$ $y_{10},$ $\dots,$ $y_{n0}$ for
the point and family of convex subsets of $P$ so obtained.
Let us also write $q_0$ for the vertex
of $x_0'$ opposite to the face $y_{10}.$

To beef these sets up to the desired members
of $\Conv(\R^n)_{\geq\{0\}},$ we now take
$$\textstyle x=\ch(x_0\cup\{0\}),\quad y_i=\ch(y_{i0}\cup
\{0\})\ \text{ for } i\neq 1,\quad y_1=\ch(y_{10}\cup\{0,q_0/2\}).$$
For the same reason as in the proof of Theorem~\ref{X.25},
the operation of intersecting with $P$ commutes with
lattice operations on these convex sets; hence
when $D\op_{n-1}$ is evaluated at
the above arguments, the left-hand side meets the plane $P$ only
in the set $x_0,$ while the right-hand side
will also contain the point $p_0.$
Now observe that $\bigwedge_i\,y_i$ will
be the line-segment $\ch(0,q_0/2).$
Letting $p=p_0,$ $q=q_0/2,$ we
see that these lie on different lines
through $0,$ and deduce from the implication (ii)(b)${\implies}$(ii)(a)
of Lemma~\ref{X.29} that $z$ and $z'$ can be chosen
so that the required inequality holds.\end{proof}

\section{Encore!}\label{NH_Radon}

Carath\'eodory's and Helly's Theorems are two members of a well-known
triad of results on convex sets in $\R^n.$
The third is\bigskip
\\
{\bf Radon's Theorem} (\cite{JE.HRC}, p.391){\bf .} {\it
Let $n \geq 0.$
Given a natural number $n,$ and $n{+}2$
points $p_1,\dots, p_{n+2}$ in $\R^n,$ there exists
a partition of $\{1,\dots,n{+}2\}$ into subsets $I_1$ and $I_2$ such
that $\ch(\{p_i\mid i\in I_1\}) \cap
\ch(\{p_i\mid i\in I_2\})\neq\emptyset.$}\qed\medskip

Can we turn this, too, into an identity for lattices of convex sets?

Yes.
First let's get rid of reference to points: Clearly an equivalent
statement is ``Given nonempty convex
sets $y_1,\dots, y_{n+2}$ in $\R^n,$ there exists
a partition of $\{1,\dots,n{+}2\}$ into
subsets $I_1$ and $I_2$ such that $(\bigvee_{i\in I_1}y_i)
\wedge (\bigvee_{i\in I_2}y_i) \neq \emptyset.\<\!$''
Next, the conclusion that such a partition exists can
be condensed into the single inequality
$\bigvee_{I_1,I_2}\,((\bigvee_{i\in I_1}y_i) \wedge
(\bigvee_{i\in I_2}y_i))\neq\emptyset,$ where the outer join is
over the $2^{n+2}-2$ partitions
of $\{1,\dots,n{+}2\}$ into two nonempty subsets.
In the proof of the theorem below,
Lemma~\ref{X.33} will be used to turn the above implication between
nonemptiness statements into a lattice relation.

In that theorem, I call an inequality $a \leq b$ that always
holds an ``identity'', since it can be rewritten as
an equation $a=a\wedge b.$
\begin{theorem}\label{X.45}
For every natural number $n,$ the identity
in $n{+}3$ variables $x, y_1,\linebreak[0] \dots, y_{n+2}$
\[\label{x.47}
\bigwedge_i\,(x\vee y_i)\ \leq\ x\,\vee\,\bigvee_{I_1,I_2}
((\bigvee_{i\in I_1}\,y_i) \wedge (\bigvee_{i\in I_2}\,y_i)),\]
where ``$\bigvee_{I_1,I_2}\!\!$'' denotes the join over all partitions
of $\{1,\dots,n{+}2\}$ into two non\-empty
subsets $I_1$ and $I_2,$ holds in
the lattice $\Conv(\R^n),$ and indeed
in\linebreak[3] $\Conv(\R^{n+1})_{\geq\{0\}},$ but not
in $\Conv(\R^{n+1}).$
\end{theorem}
\begin{proof}
We shall first prove (\ref{x.47}) in the simpler case
of $\Conv(\R^n),$ then show how to adapt the proof
to $\Conv(\R^{n+1})_{\geq\{0\}},$ and
finally, give the counterexample in $\Conv(\R^{n+1}).$

Given $x,y_1,\dots,y_{n+2}\in\Conv(\R^n)$ and
a point $p$ belonging to the left-hand side of~(\ref{x.47}),
we must show that $p$ also belongs to the right-hand side.
Let us choose $w$ as in Lemma~\ref{X.33} for the
given $x$ and $p.$
The assumption that $p$ belongs to the left-hand side of~(\ref{x.47})
says that it belongs to each of the meetands of that
expression, which by choice of $w$ means
that $w\wedge y_i$ is nonempty for all $i.$
Hence Radon's Theorem applied to those $n{+}2$ sets
says that for some $I_1,$ $I_2$ partitioning
$\{1,\dots,n{+}2\},$ we have
$\emptyset\neq (\bigvee_{I_1}w\wedge y_i) \wedge
(\bigvee_{I_2}w\wedge y_i).$
The latter set is contained in
$w\wedge(\bigvee_{I_1}y_i) \wedge(\bigvee_{I_2}y_i),$ and
the statement that this is nonempty now translates back to say
that $x\vee ((\bigvee_{I_1}y_i) \wedge (\bigvee_{I_2}y_i))$
contains $p,$ whence $p$ belongs
to the right-hand side of~(\ref{x.47}), as required.

If we are given
$x,y_1,\dots,y_{n+2}\in\Conv(\R^{n+1})_{\geq\{0\}},$ and a point $p$ on
the left-hand side of~(\ref{x.47}), we begin in the same way,
translating the hypothesis and desired conclusion to the
same statements about the sets $w\wedge y_i$
(though the convex set $w$ will not in general
belong to $\Conv(\R^{n+1})_{\geq\{0\}},$ and
so neither will these intersections).
This time we apply Radon's Theorem in $n{+}1$ dimensions to
the sets $w\wedge y_i$ together with $\{0\}.$
In the partition given by that theorem, let us assume without
loss of generality that $\{0\}$ goes into the second join; thus
we get a relation $\emptyset\neq(\bigvee_{I_1}w\wedge y_i)\,\wedge\,
(\{0\}\vee \bigvee_{I_2}w\wedge y_i).$
Here the first join is contained in $w,$ hence
so is the whole set, so that set is contained in
the meet of $w$ with the larger
set $(\bigvee_{I_1}y_i)\,\wedge\,(\{0\}\vee \bigvee_{I_2}y_i).$
Since all the $y$'s belong
to $\Conv(\R^{n+1})_{\geq\{0\}},$ the
joinands $\{0\}$ are now redundant, so we have again shown
that $w\wedge (\bigvee_{I_1}y_i) \wedge (\bigvee_{I_2}y_i)$ is
nonempty for some partition $\{1,\dots,n{+}2\}=I_1\vee I_2,$
which, as before, yields the desired conclusion by
choice of $w.$

To show that~(\ref{x.47}) does not hold
in $\Conv(\R^{n+1}),$ we begin essentially as in
the first counterexample in the proof of
Theorem~\ref{X.40}, letting $x$ be an \!$(n{+}1)$\!-simplex
in $\R^{n+1},$ $p$ a point outside
that simplex, and $x'$ the central reflection of $x$ through $p.$
This time, however, we let $y_1,\dots,y_{n+2}$ be
singletons, whose unique elements are the vertices of $x'.$
As in the earlier example we find that one side of the
identity in question (in this case the left-hand side
of~(\ref{x.47})) contains $p,$ while the other
is simply $x,$ since all joinands in the ``big join'' on
that side are empty; so the right-hand side does not
majorize the left-hand side.\end{proof}

I have not tried to fill in this picture, as I did with
the $D_n$ and $D_n\op,$ by looking for
related identities that would distinguish all
terms of (\ref{x.14}); but I expect that these exist.

Another observation, relevant to the whole development up
to this point, which I have not followed up on because it occurred to
me late in the preparation of this paper, is that if one
defines $\Conv(\R^n,\mathrm{cone})\subseteq \Conv(\R^n)_{\geq 0}$ to be
the sublattice consisting of those elements which are unions of rays
through $0,$ then $\Conv(\R^{n-1})$ embeds
naturally in $\Conv(\R^n,\mathrm{cone}),$ yielding an
extension of the chain of varieties of Lemma~\ref{X.13}.
It would be interesting to know whether for $n\,{>}\,1$ the varieties so
interpolated are distinct from those that precede and follow them.

\section{The sublattice of compact convex sets.}\label{NHcp}

From Carath\'eodory's Theorem, we see that
\[\label{x.48}
\text{For any compact subset }S \subseteq\R^n,\ \ \ch(S)
\text{ is also compact.}\]
Hence the join in $\Conv(\R^n)$ of two compact subsets
is compact, hence the set of compact convex subsets
of $\R^n$ (often called ``convex bodies'' in the
literature, e.g., \cite{HbkCnv}, \cite[\S3.1]{KA.CvG},
\cite[\S12]{FW+MS}) is a
sublattice $\Conv(\R^n,\,\cp)\subseteq \Conv(\R^n).$
An obvious question is how the identities of this sublattice
compare with those of $\Conv(\R^n);$ i.e., whether
it satisfies any identities that the larger lattice does not.
Huhn \cite[proof of Lemma~3.1]{AH.n-d} answered this question
in the negative, by showing that the still smaller lattice of polytopes
(convex hulls of finite sets) does not.
Let me give a slightly different proof of the same result.
Huhn used the fact that an intersection of polytopes is a
polytope, but the next result is applicable to a finitary closure
operator that need not have the property that an intersection of
closures of finite sets is again one.
(The meaning of ``finitary'' was recalled
in the paragraph preceding Lemma~\ref{X.24};
the notation $\V(\,...\,)$ used below was defined in~\S1.)
\begin{proposition}\label{X.50}
Let $\cl$ be a finitary closure operator on a set $X,$ and
let $L$ be the lattice of subsets of $X$ closed under $\cl.$
Then every lattice relation satisfied by all families of elements
of $L$ that are closures of finite subsets
of $X$ is an identity of $L.$
In particular, if $L'$ is any sublattice of $L$ which
contains all closures of finite subsets
of $X,$ then $\V(L')=\V(L).$
\end{proposition}
\begin{proof}
Let us topologize the power set $2^X$
by taking as a basis of open sets the sets
$$\textstyle
U(A,B)\ =\ \{Y\in 2^X\mid A\subseteq Y\subseteq B\},$$
where $A$ ranges over the finite subsets
of $X,$ and $B$ over arbitrary subsets.
This is stronger than the usual power-set topology, which
only uses the sets of the above form with $B$ cofinite.
Thus, our topology is Hausdorff, though not in general compact;
hence its restriction to $L$ is also Hausdorff, with
basis of open sets given by the sets $U_L(A,B)=U(A,B)\cap L.$
(Of course, $U_L(A,B)$ is nonempty only when $\cl(A)\subseteq B.)$

We claim that under this topology, the lattice operations of $L$ are
continuous, and the closures of finite subsets of $X$ are dense in $L.$
This will imply that any lattice identity holding on
that dense subset must hold on all of $L,$ from which
the final conclusion will clearly follow.

To see that closures of finite sets are dense, note that
every nonempty basic set $U_L(A,B)$ contains the element $\cl(A),$
which is such a set.

The continuity of the meet operation is also straightforward:
If $x,y\in L$ are such that $x\wedge y,$ i.e., $x\cap y,$ lies in
$U_L(A,B),$ then $U_L(A,x\cup B)$ and $U_L(A,y\cup B)$ are neighborhoods
of $x$ and $y$ respectively such that the
intersection of any member of the first neighborhood and
any member of the second lies in $U_L(A,B).$

Finally, suppose $x,y\in L$ are such that
$x\vee y,$ i.e., $\cl(x\cup y),$ lies in $U_L(A,B).$
Then the finite set $A$ is contained
in $\cl(x\cup y),$ so by finitariness of $\cl,$ there is a
finite subset of $x\cup y$ whose closure contains
all elements of $A;$ let us write this subset
as $A_x\cup A_y,$ where
$A_x \subseteq x$ and $A_y \subseteq y.$
Then $U_L(A_x,x)$ and $U_L(A_y,y)$ will
be neighborhoods of $x$ and $y$ respectively such that the
join of any member of the first neighborhood and
any member of the second is a member of $U_L(A,B).$
(In the power-set topology, $\vee$ is generally discontinuous;
this is why we needed a different topology.)\end{proof}

For $X=\R^n$ and $\cl=\ch,$ the $L$ of the
above proposition is $\Conv(\R^n).$
Since the convex hull of a finite
set is compact, we can apply the last sentence of the
proposition with $L'=\Conv(\R^n,\,\cp),$ getting the first statement of
the next theorem.
Taking $\cl=\ch(\{0\}\cup-)$ we similarly get the second.

\begin{theorem}[Huhn]\label{X.52}
For every natural number $n,$
$\V(\Conv(\R^n))=\linebreak[3]\V(\Conv(\R^n,\,\cp)),$ and
$\V(\Conv(\R^n)_{\geq\{0\}})=\V(\Conv(\R^n,\,\cp)_{\geq\{0\}}).$\qed
\end{theorem}

However, $\Conv(\R^n,\,\cp)$ is known also to have interesting
elementary properties not possessed by $\Conv(\R^n).$
Let us recall that an {\em extremal point} of
a convex set means a point which is not in the convex hull of any two
other points of the set, and the
following result (\cite[p.276]{HbkCnv}).\bigskip
\\
{\bf Theorem}\ \ (Minkowski){\bf .} {\it
Every compact convex subset of $\R^n$ is the
convex hull of its set of extremal points}.\qed\medskip

Recall also that a lattice $L$ is called
{\em join semidistributive} if for all $x, y_1, y_2\in L$ one has
\[\label{x.54}
x\vee y_1=x\vee y_2\ \ \implies\ \ x\vee y_1=x\vee(y_1\wedge y_2).\]

\begin{lemma}[{$=$\cite[Theorem~3.4]{KA.CvG},
generalizing \cite[Theorem~15]{B+B.po}}]\label{X.56}
For every positive integer $n,$ $\Conv(\R^n,\,\cp)$ is
join semidistributive.
\end{lemma}
\begin{proof}
Every extremal point of a join $x\vee y$ must belong to $x$ or to
$y,$ since by definition it cannot arise as a convex combination of
other points of $x\vee y,$ hence if $x\vee y_1=x\vee y_2$ as
in the hypothesis of~(\ref{x.54}), extremal points of this set that do
not belong to $x$ must belong to $y_1,$ and likewise to $y_2.$
Thus every extremal point of $x\vee y_1=x\vee y_2$ belongs
to $x\cup(y_1\cap y_2);$ hence
by Minkowski's Theorem the convex hull of the latter
set, $x\vee(y_1\wedge y_2),$ contains the former set.
The reverse inclusion is trivial.\end{proof}

On the other hand, for $n\geq 2,$ $\Conv(\R^n)$ is not
join semidistributive; indeed, the next result will show the failure of
successively weaker properties, beginning with
join-semidistributity, as $n$ increases.
Following Geyer \cite{WG.ksd}, let us say that a
lattice $L$ is {\em \!$n$\!-join semidistributive} for a
positive integer $n$ if for all $x, y_1, \dots, y_{n+1}\in L$ one has
\[\label{x.58}
x\vee y_1=\dots=x\vee y_{n+1}\ \ \implies\ \ x\vee y_1=
x\vee(\bigvee_{i,j;\ i\,\neq\,j}\ y_i\wedge y_j).\]
Thus, join semidistributivity is the $n=1$ case.
One defines {\em \!$n$\!-meet semidistributivity} dually.

The \!$(n{+}4)$\!-element lattice of
height~$2,$ $M_{n+2},$ with least element $0,$ greatest
element $1,$ and $n{+}2$ incomparable
elements $y_1,\dots,y_{n+2},$ is neither \!$n$\!-join semidistributive
nor \!$n$\!-meet semidistributive, as may be seen
by putting $y_{n+2}$ in the role of $x$ in~(\ref{x.58}),
and in the dual statement.
We shall now see that there are several sorts of sublattices with
that structure within the lattices $\Conv(\R^n).$
The ``open bounded'' case of the next result was shown to me by
D.\,Wasserman.
\begin{lemma}[Wasserman and Bergman]\label{X.60}
For every $n\,{>}\,1,$ $\Conv(\R^n)$ contains
copies of $M_{n+1}$ consisting of open bounded
sets and copies consisting of closed unbounded sets,
in both cases with least element $\emptyset;$ and
$\Conv(\R^n)_{\geq\{0\}}$ contains
copies of $M_{n+1}$ consisting of bounded sets, with
least element $\{0\}.$
Also, $\Conv(\R^2)_{\geq\{0\}}$ \r(and
hence $\Conv(\R^n)_{\geq\{0\}}$ for all $n\,{>}\,1)$ contains
copies of $M_c,$ the height-$2$ lattice of continuum
cardinality \r(and hence contains copies of its
sublattices $M_m$ for all natural numbers $m)$ consisting
of vector subspaces, with least element $\{0\}.$

In particular, for $n>1$ neither $\Conv(\R^n)$
nor $\Conv(\R^n)_{\geq\{0\}}$ is \!$N$\!-join or \!$N$\!-meet
semidistributive for any $N,$ and the sublattices of open bounded
sets in $\Conv(\R^n)$ and of bounded
sets in $\Conv(\R^n)_{\geq\{0\}}$ are
not \!$(n{-}1)$\!-join or \!$(n{-}1)$\!-meet semidistributive.
\end{lemma}
\begin{proof}
To get a copy of $M_{n+1}$ consisting of open bounded
sets, start with any \!$n$\!-simplex $x,$ let
its faces be $x_1, \dots, x_{n+1},$ choose
an interior point $p$ of $x,$ and for
each $i$ let $y_i$ be the interior
of the \!$n$\!-simplex $\{p\}\vee x_i$ (or any open convex subset of
that \!$n$\!-simplex which has the whole face $x_i$ in
its closure).
We see that the join of any two of the $y$'s will have
in its closure two faces of $x,$ hence all vertices
of $x,$ hence its closure must be $x,$ hence being
itself open and convex, it must be the interior of $x.$
On the other hand, the pairwise intersections of the $y_i$ are all
empty.
Hence the lattice generated by these sets is isomorphic to $M_{n+1}.$

For the closed unbounded example with least element
again $\emptyset,$ extend each
of the sets in the preceding example to an infinite cone with
apex $p,$ displace each of these cones away
from $p$ (say by translating it by the vector from the opposite
vertex of $x$ to $p),$ and take their closures.
Then every pairwise join is seen to be the whole
of $\R^n,$ while every pairwise meet is again empty.

For the bounded example in $\Conv(\R^n)_{\geq\{0\}},$ take the
first example above, assuming $p=0,$ and use
as our new $y_i$ the union of the $y_i$ of that example with $\{0\}.$
(Thus, $0$ will be the unique boundary-point belonging
to each of these sets.)

Finally, an $M_c$ in $\Conv(\R^2)_{\geq\{0\}}$ is given by
the set of all lines through $0.$\end{proof}

There are also cases where we can show the
failure of \!$m$\!-join
semidistributivity in a natural lattice of convex sets,
but where that lattice probably does not contain a copy of $M_{m+2}.$
To get such an example in the lattice of convex open subsets
of $\R^2$ containing $0,$ for
any $m>1,$ take $m{+}2$ distinct lines
through $0$ and ``thicken'' these to open
sets $x, y_1, \dots, y_{m+1}$ of width~$1.$
The lattice that these generate will
not be $M_{m+2},$ but clearly fails to satisfy~(\ref{x.58}).
We can get such examples for bounded open sets in $\R^2,$ though
in this case they fail to have a common
point: Fix a triangle $T,$ and let $x$ be the interior of
any triangle lying inside $T$ and sharing one edge
with $T,$ but not the opposite vertex $p.$
Then take $m{+}1$ ``small narrow'' triangles
inside $T$ that have $p$ as a common vertex but no other
point in common, and
the convex hull of whose union is disjoint from $x,$ and
let $y_1, \dots, y_{m+1}$ be their interiors.

The parenthetical comment in the first sentence of the proof of the
Lemma~\ref{X.60} shows that the shapes of the convex sets forming
a copy of $M_{n+1}$ in the lattice of open bounded convex
sets are not unique; but I don't know an example where the top
element of such a sublattice is not an open \!$n$\!-simplex.
It would also be of interest to know whether this lattice contains
copies of $M_{n+2},$ and if not,
whether it is \!$n$\!-join or \!$n$\!-meet semidistributive.
We shall obtain a few related results in subsequent sections.

J\'onsson
and Rival \cite[Lemma~2.1]{BJ+IR} show that
a lattice is join and meet semidistributive if and only if
two auxiliary overlattices contain no isomorphic copies of any
member of a certain list of 6 lattices, beginning with $M_3.$
The above ``small narrow triangle'' construction gives,
when $m=2,$ a copy of the lattice $L_4$ of their list.

Incidentally, Geyer's concept of \!$n$\!-join semidistributivity,
which we have been using, does not have any obvious relationship with
Huhn's \!$n$\!-distributivity.
Although for $n\,{=}\,1$ they give the conditions of
join-semidistributivity and distributivity respectively, of
which the latter implies the former, no such implication
holds for larger $n.$
For instance, the
lattice $M_c$ is \!$2$\!-distributive in Huhn's sense, but it
is not \!$N$\!-join semidistributive for
any natural number $N$ in Geyer's sense.

\section{Open bounded sets do not satisfy additional identities.}

Let us denote by $\Conv(\R^n,\,\ob)$ the
lattice of open bounded convex subsets of $\R^n.$
We shall show that this lattice,
like $\Conv(\R^n,\,\cp),$ satisfies the same
identities as $\Conv(\R^n).$
The idea is that compact sets can be approximated by
open bounded sets containing them, from which we shall deduce
that any identities of $\Conv(\R^n,\,\ob)$ are
also identities of $\Conv(\R^n,\,\cp),$ and so by
Theorem~\ref{X.52} are identities of $\Conv(\R^n).$

To approximate compact sets by open sets, we need a different topology
on $2^{\R^n}$ from the one used earlier;
let us again describe this in a general context.
If $X$ is any topological space, we may
topologize $2^X$ using a basis of
open sets with the same form as before,
\[\label{x.62}
U(A,B)\ =\ \{Y\in 2^X\mid A\subseteq Y\subseteq B\},\]
but where, this time, $A$ ranges over all subsets
of $X,$ while $B$ is restricted to open subsets.
If $X$ is Hausdorff (or even T1),
we see that this family of open sets again includes those defining
the power-set topology, so our topology is again Hausdorff.
Note that for each $A\in 2^X,$ the
sets~(\ref{x.62}) with $A$ as first argument and with second argument
containing $A$ form a neighborhood
basis of $A$ in $2^X.$
Thus in proving the next lemma, we shall take it as understood that
to be ``sufficiently close to'' a set $A$ means
to contain $A$ and be contained in some specified open
neighborhood $B$ of $A$ in $X.$

In the formulation of that lemma, note that the statement that a
function is continuous at arguments with a given property does not
simply mean that the restriction of the function to the set of such
arguments is continuous, but, more, that such arguments are points of
continuity of the whole function.
\begin{lemma}\label{X.63}
If $X$ is any topological space, then the binary operation
$\cup: 2^X\times 2^X \to 2^X$ is continuous in the topology described
above; if $X$ is normal \r(i.e., if disjoint closed subsets of $X$
have disjoint open neighborhoods\/\r) then the binary operation
$\cap$ is continuous at arguments given by pairs of closed sets, and
if $X=\R^n$ with the usual topology, then the unary convex-hull
operation $\ch: 2^X \to 2^X$ is continuous at compact sets.
\end{lemma}
\begin{proof}
To show continuity of $\cup,$ consider sets $A_1,$ $A_2,$
and an open neighborhood $B$ of $A_1\cup A_2$ in $X.$
Then we see that the union of any member
of $U(A_1,B)$ and any member of $U(A_2,B)$ contains $A_1\cup A_2$ and
is contained in $B,$ as required.

For the case of $\cap,$ let $A_1$ and $A_2$ be
closed sets, and $B$ any open neighborhood of $A_1\cap A_2$ in $X.$
Then $A_1-B$ and $A_2-B$ are disjoint closed sets,
hence they have disjoint open
neighborhoods $C_1$ and $C_2.$
We see that $B\cup C_1$ and $B\cup C_2$ will
be open neighborhoods of $A_1$ and $A_2$ which intersect in $B$ (by
distributivity of the lattice $2^X),$ and it follows that
the intersection of a member of $U(A_1,B\cup C_1)$ and a member
of $U(A_2,B\cup C_2)$ will belong to $U(A_1\cap A_2,B),$ as required.

For the final assertion, let $A$ be a compact
subset of $\R^n,$ and $B$ any open neighborhood of $\ch(A).$
By compactness of $\ch(A),$ there is some $\epsilon>0$ such
that the set $C$ of all points having
distance $< \epsilon$ from $\ch(A)$ is contained in $B.$
This set $C$ is a convex open neighborhood
of $\ch(A),$ hence $\ch$ carries $U(A,C)$ into
$U(\ch(A),C) \subseteq U(\ch(A),B),$ as required.\end{proof}

(I played with several topologies before getting the one that
made the above result -- in particular, continuity of
intersection -- easy to prove.
Some of these might be preferable for other considerations of
the same sort.
Under the above topology, every open set $A$ is
an isolated point, since $U(A,A)$ is a singleton.
If one wants to approximate open sets by larger open sets,
one might prefer a weaker topology in which the
conditions on $Y$ in~(\ref{x.62}) are, say, strengthened
to $A \subseteq Y,$ $\cl(Y)\subseteq B.)$

It follows from Lemma~\ref{X.63} that in the topology we have defined,
the lattice operations of $\Conv(\R^n)$ are continuous at
arguments belonging to $\Conv(\R^n,\,\cp).$
Moreover, $\Conv(\R^n,\,\cp)$ lies in the closure
of $\Conv(\R^n,\,\ob),$ since every
compact convex set $A$ is the limit in this
topology, as $\epsilon \to 0,$ of
the open convex set of points at distance $< \epsilon$ from $A.$
Hence any lattice identities
holding in $\Conv(\R^n,\,\ob)$ also
hold in $\Conv(\R^n,\,\cp).$
The same considerations apply to the pair of lattices
$\Conv(\R^n,\,\ob)_{\geq\{0\}}$ and
$\Conv(\R^n,\,\cp)_{\geq\{0\}}.$
In view of Theorem~\ref{X.52}, these observations give us
\begin{theorem}\label{X.68}
For every natural number $n,$
$\V(\Conv(\R^n))=\V(\Conv(\R^n,\,\ob)),$ and
$\V(\Conv(\R^n)_{\geq\{0\}})=\V(\Conv(\R^n,\,\ob)_{\geq\{0\}}).$\qed
\end{theorem}

Let us note that there is an order-preserving
bijection between the elements of $\Conv(\R^n,\,\ob)$ and
the compact convex subsets of $\R^n$ having
nonempty interior, given by the operation of topological closure,
with inverse given by topological interior.
This is not, however, an isomorphism between sublattices
of $\Conv(\R^n),$ because the class of compact convex
sets with nonempty interior is not closed under intersection.
(E.g., consider two adjacent closed polygons in $\R^2.$
Nor can we get around this problem by going to a
homomorphic image of $\Conv(\R^n,\,\cp)$ where
sets without interior are identified with $\emptyset,$ since the join
of two such sets can have nonempty interior.)

On the other hand, the set of compact convex sets which are
neighborhoods of $0$ {\em is} a sublattice
of $\Conv(\R^n,\,\cp)_{\geq\{0\}},$ and the above
correspondence gives us an isomorphism between it
and $\Conv(\R^n,\,\ob)_{\geq\{0\}};$ so we can apply results about
$\Conv(\R^n,\,\cp)_{\geq\{0\}}$ to $\Conv(\R^n,\,\ob)_{\geq\{0\}}.$
(More generally, any sublattice of $\Conv(\R^{n},\,\ob)$ whose members
have a common point $p$ is contained in
$\Conv(\R^n,\,\ob)_{\geq\{p\}} \cong\Conv(\R^n,\,\ob)_{\geq\{0\}},$
and so can be studied in the same fashion.)
Thus, despite the examples of Lemma~\ref{X.60}, we have
\begin{corollary}[to Lemma~\ref{X.56}]\label{X.69}
$\Conv(\R^n,\,\ob)_{\geq\{0\}}$ is join semidistributive.
Hence, every counterexample to join semidistributivity
in $\Conv(\R^n,\,\ob)$ has the property
that the intersection of the three sets involved is empty.\qed
\end{corollary}

Since we are considering elementary properties in which $\Conv(\R^n)$
and\linebreak[3] $\Conv(\R^n,\,\ob)$ agree or differ, we should note
the obvious difference, that the former is atomistic (every element
is a possibly infinite join of atoms), while the latter has no atoms.
Cf.\ \cite{MB.bia} and papers referred to there, in which
lattices of convex sets and related structures are characterized in
terms of properties of their atoms, and also \cite{KA.CvG}.

We noted earlier that lattices $\Conv(\R^n)_{\geq S}$ for
nonempty $S$ satisfy all identities holding
in $\Conv(\R^n)_{\geq\{0\}}.$
Let us end this section by using Theorem~\ref{X.68} to show that for
bounded $S,$ the converse is also true.
\begin{corollary}[to Theorem~\ref{X.68} and
proof of Theorem~\ref{X.52}]\label{X.70}
For every natural number $n$ and every bounded set $S \subseteq \R^n,$
$$\textstyle
\V(\Conv(\R^n)_{\geq S})\ =\ \V(\Conv(\R^n,\ \cp)_{\geq S})$$
$$\textstyle
=\ \V(\Conv(\R^n,\ \ob)_{\geq S})\ =\ \V(\Conv(\R^n)_{\geq\{0\}}).$$
\end{corollary}
\begin{proof}
By a translation, we can assume without loss of generality
that $0\in S;$ thus the last of the above lattices
contains all the others, so letting
$f(x_1,\dots,x_m)=g(x_1,\dots,x_m)$ be any identity
not satisfied there, it suffices to prove that it
is not satisfied in any of the other lattices.

Now Theorem~\ref{X.68} shows that $f=g$ is not an identity
of $\Conv(\R^n,\,\ob)_{\geq\{0\}},$ hence
we can choose open bounded convex
sets $x_1,\dots,x_m$ containing $0$ which do not satisfy it.
The intersection of these sets is a neighborhood of the
origin, and dilating the $x_i$ by a large enough real
constant, we can assume without loss of generality
that this neighborhood contains $S.$
Hence $f=g$ is also not an identity
of $\Conv(\R^n,\,\ob)_{\geq S},$ hence
not an identity of $\Conv(\R^n)_{\geq S}$ either.

The case of $\Conv(\R^n,\,\cp)_{\geq S}$ is similar.
Again take
$x_1,\dots,x_m\in\linebreak[3]\Conv(\R^n,\,\ob)_{\geq\{0\}}$ not
satisfying $f=g.$
We saw in the proof of Theorem~\ref{X.52} that if we
approximate $x_1, \dots, x_m$ closely
enough from below in the topology of that proof by
compact convex subsets $y_1, \dots, y_m,$ these
approximating sets also fail to satisfy that identity.
Since the intersection of the $x_i$ is a neighborhood of $0,$ it
contains an \!$n$\!-simplex with $0$ in its interior;
so we can take all the $y_i$ to contain the finitely many vertices
of that simplex, hence to be neighborhoods of $0.$
As before, we may now dilate them so that they all
contain $S,$ getting the required result.\end{proof}

\section{The possibility of surface phenomena.}

The technique by which we just proved Corollary~\ref{X.70} can be
inverted to show that if $T \subseteq \R^n$ is any convex
set with nonempty interior, then the lattice of convex
sets contained in $T,$ and its sublattices of
compact convex subsets of $T$ and open bounded convex
subsets of $T,$ satisfy the same identities as $\Conv(\R^n).$
Namely, given any identity
not holding in $\Conv(\R^n),$ we already know that
we can find elements $x_1,\dots,x_m$ not
satisfying it in $\Conv(\R^n,\,\cp).$
Assuming without loss of generality that $0$ lies in the
interior of $T,$ we can shrink $x_1,\dots,x_m$ by a constant,
so that they lie in $T$ as well.
A similar argument shows that the lattice of
convex subsets of $T$ which contain a specified
point of the interior of $T$ satisfies the same
identities as $\Conv(\R^n)_{\geq\{0\}}.$

However, if we specify two convex sets $S \subseteq T,$ say
with $S$ compact and $T$ open, and look at the
interval $[S,\ T]=\{x\in\Conv(\R^n)\mid S \subseteq x \subseteq T\},$ it
is not clear whether, for some choices
of $S$ and $T,$ this may satisfy more identities
than hold in $\Conv(\R^n)_{\geq\{0\}}.$
(Picturing $S$ and $T$ as ``very close'', e.g.,
the closed ball of radius $1$ and the open ball of
radius $1+\epsilon,$ explains, I hope, the title of this section.)

Let us relax our assumptions on $S$ and $T$ for
a moment and look at a more extreme example.
If we take for $S$ the open
unit ball and for $T$ the closed unit ball in $\R^n,$ then
every set $x$ with $S \subseteq x \subseteq T$ is convex,
so in this case $[S,\ T]$ may be identified with the lattice
of all subsets of the unit sphere, which is distributive, though we
have seen that $\Conv(\R^n)_{\geq\{0\}}$ is not
even \!$(n{-}1)$\!-distributive.

In the case where $S$ is compact and $T$ open,
however, things cannot go that far:
\begin{lemma}\label{X.80}
Let $S \subseteq T$ be convex subsets of $\R^n.$
If $S$ is compact and $T$ is open and nonempty, or more generally, if
some hyperplane $P \subseteq \R^n$ disjoint from $S$ intersects $T$ in
a set with nonempty relative interior \r(i.e., is such
that $P\cap T$ is \!$(n{-}1)$\!-dimensional\/\r),
then $\V([S,T]) \supseteq\V(\Conv(\R^{n-1})).$
\end{lemma}
\begin{proof}
Clearly, a pair $S \subseteq T$ satisfying the first hypothesis
satisfies the second, so let us assume the latter.
Let $H$ be the closed half-space of $\R^n$ bounded
by $P$ that contains $S$ (or if $S$ is empty, either
of the closed half-spaces bounded by $P),$ and consider
the sublattice $[S,\,H\cap T] \subseteq [S,T].$
The operation of intersecting with $P$ can be seen to give a lattice
homomorphism $[S,\,H\cap T] \to [\emptyset,P\cap T],$ and this is
surjective, since it has the set-theoretic section $x \mapsto S\vee x.$
Since $P\cap T$ has nonempty interior in $P,$ the observation
in the first paragraph of this section shows that $[\emptyset,P\cap T]$
satisfies precisely the identities of $\Conv(\R^{n-1}).$
Hence $[S,\,H\cap T],$ since it maps homomorphically
onto $[\emptyset,P\cap T],$ cannot satisfy any identities not
satisfied by $\Conv(\R^{n-1}),$ so neither
can the larger lattice $[S,T].$\end{proof}

It is not evident whether, for $S$ compact and $T$ open, $\V([S,T])$
can ever be strictly smaller than $\V(\Conv(\R^n)_{\geq\{0\}});$
nor, for that matter, whether it can ever fail to
be strictly smaller, if $S$ has nonempty interior and $T$ is bounded.
Another interval $[S,T]\subseteq\Conv(\R^n)$ whose identities it
would be interesting to investigate is given by letting
$S=\{0\}$ and $T$ be a closed half-space with $0$ on its boundary.
Again these identities must lie somewhere between those
of $\Conv(\R^{n-1})$ and $\Conv(\R^n)_{\geq\{0\}}.$

In this and preceding sections we have used from time to
time the fact that translations and dilations preserve convexity.
More generally, if $\phi$ is any projective transformation on
\!$n$\!-dimensional projective space $\mathbb P^n \supseteq \R^n,$
then convex subsets of $\R^n$ which do not meet the hyperplane
that $\phi$ sends to infinity are taken by $\phi$ to convex sets.
Hence if $S\subseteq T$ are such convex sets,
$\phi$ induces a lattice isomorphism $[S,T]\cong [\phi(S),\phi(T)].$
This observation might be useful in classifying
varieties generated by such intervals.

Another sort of sublattice of $\Conv(\R^n)$ that it
might be interesting to investigate is that of all convex sets
that are carried into themselves by a given affine map; e.g.,
the orthogonal projection onto a specified subspace.
(If that subspace is $\{0\},$ we get $\Conv(\R^n)_{\geq\{0\}}$
with one additional element~$\emptyset$ thrown in.)

\section{Dualities.}\label{NHDual}

Let us recall the definition of a concept we have
referred to a couple of times in passing.
A {\em closed half-space} in $\R^n$ means
a set of the form
\[\label{x.71}
\{p\in\R^n\mid f(p)\,\leq\,\lambda\},\]
for some nonzero linear functional $f$ on $\R^n$ and
some real number $\lambda.$
It is a standard result that every closed convex subset
of $\R^n$ is an intersection of closed half-spaces.

The half-space (\ref{x.71}) contains the point $0$ if and only
if $\lambda$ is nonnegative, hence closed convex sets containing $0$ can
be characterized as intersections of half-spaces (\ref{x.71})
having $\lambda \geq 0.$
Such sets can, in fact, be expressed as intersections of such
half-spaces with $\lambda>0,$ since a half-space~(\ref{x.71})
with $\lambda=0$ is the intersection of
all the half-spaces with the same $f$ and positive $\lambda.$
But a half-space~(\ref{x.71}) with $\lambda$ positive can be written
as $\{p\in\R^n\mid \lambda^{-1}f(p) \leq 1\},$ or, expressing the
linear functional $\lambda^{-1}f$ as the dot product with
some $q\in\R^n,$ as
\[\label{x.72}
\{p\in\R^n\mid q\<{\cdot}\<p\,\leq\,1\}.\]
Thus, for any subset $S \subseteq \R^n,$ if we define
\[\label{x.74}
S^*\ =\ \{q\in\R^n\mid (\forall p\in S)\ q\<{\cdot}\<p\,\leq\,1\},\]
then $S^{**}$ will be the least closed convex
set containing $S\cup\{0\}.$
Moreover, we see that $S^*$ will also be a closed
convex subset containing $0,$ which uniquely determines and
is determined by $S^{**}.$
Thus, the operator $^*$ gives a bijection of the family of all closed
convex sets containing $0$ with itself, which is easily seen to
be inclusion-reversing.
(This is an example of a Galois connection;
cf.\ \cite[\S5.5]{GB.245} for a general development
of the concept, with many examples.)
Let us call two closed convex sets containing $0$ that are
related in this way dual to one another.
(The dual of a convex set is sometimes called its polar set,
e.g., in \cite{JL.Clu}.)
Examples in $\R^3$ are a cube and an octahedron of appropriate radii
centered at the origin, and similarly a dodecahedron and an icosahedron.
The unit sphere is self-dual.

The class of closed convex subspaces of $\R^n$ forms a lattice
(by general properties of Galois connections), which, like
$\Conv(\R^n)_{\geq\{0\}},$ has intersection as its meet operation; but
the join operations do not everywhere coincide -- a consequence of the
fact that, given closed sets $u,v\in\Conv(\R^n)_{\geq\{0\}},$
their join in that lattice, $\ch(u\cup v),$ may not be closed,
so that to get their join as closed convex sets,
one must take its topological closure.
For example, let $n=2,$ and let $u$ be the closed strip
$\{(x,y)\mid -1{\,\leq\,}y{\,\leq\,}1\}$ and $v$ the line
segment $\ch(0,\,(0,2)).$
Then the join of $u$ and $v$ in $\Conv(\R^n)_{\geq 0}$ is
$$\{(x,y)\mid -1\leq y<2\}\ \,\cup\ \,\{(0,2)\},$$
while their join in the corresponding lattice of closed convex sets
is $\{(x,y)\mid -1\,{\leq\,}\,y\,{\leq}\,2\}.$
In view of this difference in operations, care is needed when
using our duality on closed convex sets to deduce results about
the lattice $\Conv(\R^n)_{\geq\{0\}}.$

If $x$ and $y$ are mutually dual closed convex sets
containing $0$ in $\R^n,$ it is not hard to see
that one of them is bounded (i.e., compact) if and only
if $0$ is an interior point of the other.
It follows that the class of closed bounded convex sets having $0$ in
their interior is self-dual; moreover, we saw
at the beginning of~\S\ref{NHcp} that the join
in $\Conv(\R^n)$ of two compact sets is again compact, from
which it follows that unlike the lattice of all closed convex
sets containing $0,$ this is a sublattice of $\Conv(\R^n).$
It is easy to see from the equality of the second
and fourth varieties in Corollary~\ref{X.70} that this sublattice
satisfies the same identities as $\Conv(\R^n).$
Hence the existence of the anti-automorphism just noted gives us

\begin{theorem}\label{X.73}
The class of lattice identities satisfied
by $\Conv(\R^n)_{\geq\{0\}}$ is self-dual, i.e., closed
under interchanging all instances of $\vee$ and $\wedge.$
Equivalently, the variety $\V(\Conv(\R^n)_{\geq\{0\}})$ is closed
under taking dual lattices.\qed
\end{theorem}

And indeed, the identities proved
for $\Conv(\R^n)_{\geq\{0\}}$ in
Theorems~\ref{X.25} and~\ref{X.40} respectively are dual to one another.
This is not true of the identities proved
for $\Conv(\R^n)$ in those theorems, as is easily verified:
\begin{exercise}\label{X.77}
{\rm(i)} For every positive integer $n,$ show by example
that $\Conv(\R^n)$ does not satisfy either of the identities
$$\textstyle((D_n\op)\wedge z)\vee y_1\vee y_2\quad\text{and}\quad
(((D_n)\vee z')\wedge z)\vee(\bigvee_i\,y_i),$$
dual to the identities proved for that lattice in
Theorems~\ref{X.25} and~\ref{X.40} respectively.

{\rm(ii)} Show that these observations together with Theorem~\ref{X.73}
yield an alternative way of verifying that in each
of Theorems~\ref{X.25} and~\ref{X.40}, the identities proved for
$\Conv(\R^n)$ are not satisfied by $\Conv(\R^{n+1})_{\geq\{0\}}.$
\end{exercise}

Theorem~\ref{X.73} also implies
that $\Conv(\R^{n+1})_{\geq\{0\}},$ and hence
also $\Conv(\R^{n}),$ satisfies the dual of the identity of
Theorem~\ref{X.45}, which we had not previously obtained.\smallskip

Just as every closed convex set is an intersection of closed
half-spaces (\ref{x.71}), so every open convex set is an
intersection of open half-spaces,
$$\{p\in\R^n\mid f(p) < \lambda\}.$$
But here the converse is not true.
Indeed, every {\em closed} half-space is also an intersection of open
half-spaces, so the class of intersections of open half-spaces
includes both the open and the closed convex sets.
If for every subset $S \subseteq \R^n$ one defines
\[\label{x.75}
S\odual\ =\ \{q\in\R^n\mid (\forall p\in S)\ q\<{\cdot}\<p<1\},\]
one gets a duality theory for the class of $\odual\odual$\!-invariant
sets, i.e., sets satisfying $x\odual\odual=x,$ a larger class than
that covered by the preceding duality.
It is not hard to see that a necessary and sufficient condition for a
set $x\in\Conv(\R^{n})_{\geq\{0\}}$ to be ${\odual\odual}$\!-invariant
is that for every point $p$ of the boundary
of $x$ which does not belong to $x,$ there exist
a supporting hyperplane of $x$ through $p$ containing no point of $x.$
From this we can see that the union of the
open unit ball $B$ in $\R^n$ with any subset of its
boundary is $\odual\odual$\!-invariant; its $\odual$\!-dual is the
union of $B$ with the complementary subset of the boundary.
On the other hand, we find that the union of an open polygonal
neighborhood $P$ of $0$ in the plane with a nonempty finite (or
countable) subset $S$ of its boundary is never
$\odual\odual$\!-invariant, since by the above characterization of
$\odual\odual$\!-invariant sets, each point $q\in S$ forces the
$\odual\odual$\!-closure of $P\cup S$ to contain all points $p$ of
the open edge(s) of our polygon containing or adjacent to~$q.$

Like the $^{**}\!$-invariant sets (closed convex sets containing $0),$
the $\odual\odual$\!-invariant subsets of $\R^n$ form a
lattice by general properties of Galois connections, which has the
same meet operation as $\Conv(\R^n)_{\geq\{0\}},$ but
different join operation.
For instance, if $P$ is, as above, an open polygonal neighborhood
of $0$ in $\R^2,$ $q$ a point of its boundary,
and $Q=\ch(0,q),$ then writing $\vee$ for
the join operation of $\Conv(\R^2)_{\geq\{0\}},$ we find
that $P\vee Q=P\cup Q=P\cup\{q\},$ but as
we saw above, this set is not $\odual\odual$\!-invariant.
Rather, $(P\cup Q)\odual\odual,$ the
join of $P$ and $Q$ in the lattice of $\odual\odual$\!-invariant
sets, is obtained by attaching to $P\cup\{q\}$ the open edge(s)
containing or adjacent to $q.$
Further, the $^{**}\!$-invariant subsets of $\R^n$ do not even form a
sublattice of the the $\odual\odual$\!-invariant sets.
To see this in $\R^2,$ let
$$\textstyle u=\{(x,\,y)\mid\,|x|<1,\ y\geq\ \ 1/(1-x^2)-2\},$$
$$\textstyle v=\{(x,\,y)\mid\,|x|<1,\  y\leq -1/(1-x^2)+2\}.$$
We see that the join of $u$ and $v$ in $\Conv(\R^2)_{\geq\{0\}}$ is
$\{(x,\,y)\mid |x|<1\},$ which is open, hence
is $\odual\odual\!$-invariant, hence is
their join in the lattice of $\odual\odual\!$-invariant sets;
but it is not closed, hence is not their join in the lattice
of $^{**}\!$-invariant sets.

But again, one can apply $\odual\odual$\!-duality to sublattices
of $\Conv(\R^n)_{\geq\{0\}}$ which consist of
$\odual\odual$\!-invariant elements, since for these
the two lattice structures in question must agree.
In particular, one can verify that $\odual$\!-duality
interchanges compact and open sets, giving an anti-isomorphism
between the sublattices $\Conv(\R^n,\,\cp)_{\geq\{0\}}$ and
$\Conv(\R^n,\mathrm{open})_{\geq\{0\}}.$
(This anti-isomorphism can also be obtained by composing the
anti-isomorphism ``$^*\!$'' between $\Conv(\R^n,\,\cp)_{\geq\{0\}}$ and
the lattice of closed convex sets that are neighborhoods of the
origin, with the isomorphism noted earlier between that
lattice and $\Conv(\R^n,\mathrm{open})_{\geq\{0\}}.)$
The existence of this anti-isomorphism gives
\begin{corollary}[to Lemma~\ref{X.56}]\label{X.76}
For every positive integer $n,$
$\Conv(\R^n,\,\mathrm{open})_{\geq\{0\}}$ is meet semidistributive,
that is, satisfies the dual of {\rm(\ref{x.54})}.\qed
\end{corollary}

This explains the fact that in Lemma~\ref{X.60} and the discussion
that followed, though we saw that
$\Conv(\R^n,\,\mathrm{open})_{\geq\{0\}}$ was not join
semidistributive, we found no copies of $M_3$ in it.

Meet-semidistributivity does {\em not} hold in any of the lattices of
convex sets we have considered that are not defined so as to make
those sets all have some point (such as $0)$ in common;
for in each of these lattices, it is easy to get examples of
a set $x$ having empty intersection with each of two
sets $y_1$ and $y_2,$ but nonempty intersection with their join.
It also does not hold in $\Conv(\R^n,\,\cp)_{\geq\{0\}}$
for $n>1,$ since that lattice embeds in in $\Conv(\R^{n-1},\,\cp)$
(cf.\ proof of Lemma~\ref{X.13}), to which the above observation
applies.

While on the topic of join- and meet-semidistributivity, I will
note some questions and examples which are easier to state
now that we have named several sublattices of $\Conv(\R^n).$
Kira Adaricheva
(personal communication) has posed several questions of the following
form: If we take one of the sublattices
of $\Conv(\R^n)$ that we know to be join- or
meet-semidistributive, and extend it by adjoining,
within $\Conv(\R^n),$ a single ``nice'' outside element,
is the property in question already lost, and if so, does the resulting
lattice in fact contain copies of $M_3?$

In many cases this does happen.
For instance, the sublattice $L \subseteq\Conv(\R^2)$ generated by the
join-semidistributive sublattice $\Conv(\R^2,\,\cp),$ and the open disk,
which we shall denote $B,$ contains a copy of $M_3.$
To describe it, let $p_0,$ $p_1,$ $p_2$ be any three distinct points
on the unit circle, and for $i=0,1,2$ define the
point $q_i=1/5p_i + 2/5p_{i+1} + 2/5p_{i+2}$ (subscripts
evaluated\ \ mod~$3;$ I suggest making a sketch).
Let $x_i=(\ch(q_i,p_{i+1})\wedge B) \vee (\ch(q_i,p_{i+2})\wedge B).$
Then $x_1,\ x_2,\ x_3$ can be seen to generate a copy
of $M_3$ in $\Conv(\R^{2}),$ resembling,
though not identical to, the $n\,{=}\,2$ case
of the first example described in Lemma~\ref{X.60}.
An analogous construction gives a family of elements
of $L$ resembling the example
described immediately following that lemma, showing that $L$ is
not \!$m$\!-join semidistributive for any $m.$

As another example let $L$ be the sublattice
of $\Conv(\R^3)_{\geq\{0\}}$ generated
by $\Conv(\R^3,\,\ob)_{\geq\{0\}}$ and the cube $C=[-1,+1]^3.$
Letting $B$ denote the interior of $C,$ I
claim that for any triangle $T$ drawn on a
face $F$ of $C,$ $L$ contains the union
of $B$ with the interior of $T$ relative to $F.$
Indeed, we can find in $\Conv(\R^3,\,\ob)_{\geq\{0\}}$ an
open pyramid $P$ (with apex near $0)$
meeting no face of $C$ except $F,$ and
meeting $F$ in precisely the relative interior of $T.$
Then $(P\wedge C)\vee B$ will be the desired set.
If we construct three sets of this sort using as our $T$'s three
triangles in a common face $F $ of $C,$ whose relative
interiors form a copy of the $n\,{=}\,2$ case of the first example of
Lemma~\ref{X.60}, then the sublattice generated
by the three resulting sets $(P\wedge C)\vee B$ will
have the form $M_3$ (with $B$ as its least element;
cf.\ the last sentence of Corollary~\ref{X.69}).
The same trick can be used to get examples like those of the paragraph
following Lemma~\ref{X.60}, showing that $L$ is not
\!$m$\!-join semidistributive for any positive $m.$

On the other hand, I do not know what can be said about the lattices
gotten by
adjoining to $\Conv(\R^n,\,\ob)_{\geq\{0\}}$ a
closed ball $D$ about $0,$ respectively an open ball $B$ not
containing $0,$ except that the second of these lattices is not
meet-semidistributive, since we can find three ``fingers''
in $\Conv(\R^n,\,\ob)_{\geq\{0\}}$ which
meet $B$ in disjoint subsets $x,$ $y_1,$ $y_2,$ such
that $x\wedge(y_1\vee y_2)$ is nonempty.

\section{Relatively convex sets.}\label{NHRel}

Let $S$ be a subset of $\R^n,$ in general non-convex.
We shall call a subset $x \subseteq S$ convex {\em relative
to} $S$ if $x=\ch(x)\cap S;$ equivalently, if $x$ is the
intersection of $S$ with some convex subset of $\R^n.$
Such sets $x$ will form a
lattice $\RelConv(S),$ with meet operation given, as
in (\ref{x.04}), by intersection, but join now given by
\[\label{x.83}
x\vee y\ =\ \ch(x\cup y)\,\cap\,S.\]
Note that the map $x \mapsto \ch(x)$ gives
a bijection between the relatively convex subsets of $S,$ and
the subsets of $\R^n$ which are convex hulls of subsets
of $S,$ the inverse map being given by $- \cap S.$
Convex sets of the latter sort form a lattice with
join as in~(\ref{x.04}), but with meet operation
\[\label{x.85}
x\wedge y\ =\ \ch(x\cap y\cap S).\]
We observe that Carath\'eodory's Theorem (with the
``refinement'' given in the final sentence), regarded as a statement
about the closure operator $\ch(-)$ on $\R^n,$ entails
the same properties for the closure operator $\ch(-)\cap S$ on $S.$
Hence the proofs of Lemma~\ref{X.23} and of the positive assertions
of Theorem~\ref{X.25} immediately yield
\begin{proposition}[{cf.\ Huhn \cite[Lemma~3.2]{AH.n-d}}]\label{X.88}
If $n$ is a natural number and $S$ a subset of $\R^n,$ then
$\RelConv(S)$ satisfies the
identity $((D_n)\vee z)\wedge y_1\wedge y_2,$ and
if $p$ is a point of $S,$ $\RelConv(S)_{\geq\{p\}}$ satisfies the
identity $D_n.$\qed
\end{proposition}

However, the next lemma shows that the corresponding results fail
badly for the identities involving $D_n\op$ obtained
in Theorem~\ref{X.40}.
In this lemma, the second assertion embraces the first (plus two obvious
intermediate results not stated); however I include the first assertion
because both its statement and the example proving
it are more transparent than for the second.
\begin{lemma}\label{X.90}
For $S$ a subset of $\R^2,$ $\RelConv(S)$ need
not satisfy the identity $D_n\op$ for any positive integer $n.$
In fact, for $p$ an element of such
an $S,$ $\RelConv(S)_{\geq\{p\}}$ need not satisfy the
identity $(((D_n\op)\wedge z')\vee z)\wedge(\bigwedge_i\,y_i).$
\end{lemma}
\begin{proof}
To get the first assertion, let $S$ be the set consisting
of the unit circle and its center, $0.$
Given $n,$ let $q_1,\dots,q_{2n+2}$ be the successive
vertices of a regular \!$(2n{+}2)$\!-gon on that circle;
for $i=1,\dots,n{+}1,$ let $y_i=\{q_1,\dots,q_{n+1}\}-\{q_i\},$ and
let $x=\{q_{n+2},\dots,q_{2n+2}\}.$
Observe that each intersection of $n$ of the $y$'s consists
of one of the points $q_i,$ and that $x$ contains
the antipodal point; hence the join of $x$ with each intersection
of $n$ $y$'s contains~$0.$
On the other hand, the intersection of all the $y$'s is
empty, hence its join with $x$ is $x,$ which does not contain~$0.$
So $D_n\op$ fails for this choice of arguments.

To get an example where all the given sets have a common
element $p,$ and where, moreover, the indicated weaker
identity fails, let us
take for $S$ the same set as in the above example, with the
addition of one arbitrary point $r$ outside the circle,
at distance $\geq 2$ from $0.$
This time, let $q_1,\dots,q_{2n+6}$ be the successive
vertices of a regular \!$(2n{+}6)$\!-gon on $S,$ placed so
that $q_{n+2}$ lies on
the line connecting $0$ with the external point $r.$
For $i=1,\dots,n{+}1,$ let $y_i=\{q_1,\dots,q_{n+3}\}-\{q_i\},$
and let $x=\{q_{n+3},\dots,q_{2n+4}\}.$
We note that the set of subscripts of the $q$'s occurring in each
of these sets lies in an interval of length $<n{+}3,$ so
the absence of the point $0$ does not contradict convexity of these
sets relative to $S.$
Taking $p=q_{n+3}$ we see that all of these
sets belong to $\RelConv(S)_{\geq\{p\}}.$

As in the previous example, the intersection of
any $n$ of the $y$'s contains one of $q_1,\dots,q_{n+1},$ and
$x$ contains its antipodal point,
so that the join of $x$ with that intersection
contains~$0,$ hence so does the intersection of all these joins,
i.e., the value of the right-hand side
of $D_n\op$ at these arguments.
On the other hand, the intersection of all the $y$'s
is $\{q_{n+2}, q_{n+3}\},$ and
the union of this set with $x$ still has all subscripts lying in
an interval of length $<n{+}3,$ so the join of those
two sets, the left-hand side of $D_n\op,$ does not contain $0.$
Hence if we take $z'=\{q_{n+3},0\},$ the
intersection of the right-hand side of $D_n\op$
with $z'$ is $z'=\{q_{n+3},0\},$ while intersection of the
left-hand side with $z'$ is $\{q_{n+3}\}.$
If, finally, we let $z=\{q_{n+3},r\}$ and
take the joins of this element with those two intersections, we see
that in the first case the resulting set contains the
point $q_{n+2},$ since we assumed this to lie on the line-segment
from $0$ to $r,$ while in the second, it does not.
Since $q_{n+2}\in\bigwedge_i y_i,$ the two sets remain distinct on
intersecting with $\bigwedge_i y_i,$ showing the failure
of $(((D_n\op)\wedge z') \vee z)\wedge(\bigwedge_i\,y_i).$\end{proof}

The above example in $\R^2$ is also an example
in $\R^n$ for any $n \geq 2.$
For completeness, we should also consider dimensions $n=0$ and~$1.$
If we look at the chain of varieties corresponding
to that of Lemma~\ref{X.13}, but with
each $\V(\Conv(\R^n))$ replaced by
the variety generated by all the lattices $\RelConv(S)$ for
subsets $S \subseteq \R^n,$ and
each $\V(\Conv(\R^n)_{\geq\{0\}})$ by the variety generated by
lattices $\RelConv(S)_{\geq\{p\}},$ then
we see that the first term of this chain is still
the trivial variety and the next two still satisfy
the distributive identity (in the last case, by the final
assertion of Proposition~\ref{X.88} for $n=1).$
Since the distributive identity implies the identities of every
nontrivial variety, we conclude that allowing lattices of
relatively convex sets
has not enlarged the varieties we get at these three steps.
We have just shown the contrary from the fifth step on; this
leaves only the fourth step, i.e., the relation
between $\V(\Conv(\R^1))$ and the
variety generated by all lattices of
form $\RelConv(S)$ with $S\subseteq\R^1.$
Here again, it turns out that we have equality.
This follows from our observation
that $\RelConv(S)$ is isomorphic to the lattice
of convex hulls of subsets of $S,$ together with
\begin{lemma}\label{X.91}
Let $S$ be a subset of $\R^1.$
Then the convex hulls of subsets of $S$ form a
sublattice of $\Conv(\R^1).$
\end{lemma}
\begin{proof}
We have noted, for arbitrary $n$ and $S \subseteq \R^n,$ that
the convex hulls of subsets of $S$ are
closed under the join operation of the full lattices of convex sets,
so we need only show them closed under meets, i.e., intersections,
when $n=1.$
This comes down to showing that if $p\in\ch(x)\cap \ch(y),$
where $x$ and $y$ are relatively convex subsets
of $S,$ then $p\in\ch(x\cap y).$
The case where $p\in S$ is immediate, so assume the contrary.
The fact that $p\in\ch(x)$ then means
that $q_1<p<q_2$ for some $q_1,q_2\in x;$ similarly,
$q_3<p<q_4$ for some $q_3, q_4\in y.$
From the relative convexity of $x$ and $y$ and the
order-relations of these elements, we now see
that $\max(q_1,q_3)\in x\cap y$
and $\min(q_2,q_4)\in x\cap y.$
Hence $p\in\ch(\max(q_1,q_3),\,\min(q_2,q_4)) \subseteq
\ch(x\cap y).$\end{proof}

For further results on $\V(\Conv(\R^1))$ and
its subvarieties, see \cite{FW+MS2}.

Incidentally, the lattice $\RelConv(S),$ for $S$ the set used in first
part of the proof of Lemma~\ref{X.90} (consisting of the unit circle
and its center) shows that none of the identities
$(((D_m\op)\wedge z') \vee z)\wedge(\bigwedge_i\,y_i)$ implies
any of the identities $D_n\op.$
Indeed, we showed that $\RelConv(S)$
satisfies none of the latter identities; let us now show
that (unlike the lattice considered in the second part of
that proof) it satisfies all of the former.
It suffices to verify the strongest of these, the case $n=1.$
Writing $S=C \cup \{0\},$ where $C$ is the unit circle, we see that the
map $-\cap C: \RelConv(S) \to \RelConv(C)=2^C$ is a homomorphism; so if
two lattice expressions in $x,$ $y_1,$ $y_2,$ $z,$ $z'$ are
identically equal in distributive lattices, their values
in $\RelConv(S)$ will always agree except,
perhaps, as to whether they contain~$0.$
We now consider separately the
cases $0\in y_1\wedge y_2$ and $0\notin y_1\wedge y_2.$
In the former case, the two sides
of $D_1\op(x,y_1,y_2)$ agree in containing $0,$ hence are equal, so
a fortiori the two sides of
$(((D_1\op)\wedge z')\vee z)\wedge(y_1\wedge y_2)$ are equal.
In the latter case, neither side of the latter relation can
contain~$0,$ hence again they are equal.

\section{The snowflake.}\label{NHsnowflake}

Let us look at a particularly neat example of a lattice of
relatively convex subsets of a set $S.$

Let $p_1,\ p_2,\ p_3, -p_1, -p_2, -p_3$ be the successive
vertices of a regular hexagon in $\R^2$ centered at $0.$
Let
$$\textstyle S_1=\ch(p_1,-p_1),\quad
S_2=\ch(p_2,-p_2),\quad S_3=\ch(p_3,-p_3),$$
$$\textstyle S\ =\ S_1\,\cup\,S_2\,\cup\,S_3,$$
and let $L$ be the sublattice
of $\RelConv(S)_{\geq\{0\}}$ generated by the three
line-segments $S_1,$ $S_2$ and $S_3.$
(In view of the form of $S,$ I think of this example as
``the snowflake''.)

Every element of $L$ will clearly be centrally symmetric and
topologically closed; hence every such element has the form
\[\label{x.92}
\lambda_1\,S_1\ \cup\ \lambda_2\,S_2\ \cup\ \lambda_3\,S_3\quad
(\lambda_1, \lambda_2, \lambda_3\in[0,1]).\]
The join of $S_1$ and $S_2$ in this lattice must
have the form $S_1\cup S_2\cup\lambda\,S_3$ for some $\lambda\in(0,1).$
(Elementary geometry shows that $\lambda=1/2;$ but we don't need to
know this now, and will get it from a general formula soon.)
Intersecting this join with $S_3,$ we see that $\lambda\,S_3\in L.$
By symmetry we also have $\lambda\,S_1,\,\lambda\,S_2\in L,$ and
we see that these together generate a proper
sublattice of $L$ isomorphic to the whole lattice.
In particular, $L$ has infinite descending chains of elements,
e.g., $S_1 > \lambda\,S_1 > \lambda^2 S_1 > \dots.$

Let's figure out how to calculate in $L.$
The first thing we should find are the conditions
on $\lambda_1,$ $\lambda_2$ and $\lambda_3$ for
a set~(\ref{x.92}) to be relatively convex.
Calculation shows that for $\lambda_1,$ $\lambda_2$ and $\lambda_3$
nonzero, the points $\lambda_1\,p_1,$ $\lambda_2\,p_2$ and
$\lambda_3\,p_3$ are collinear if and only if
$\lambda_2^{-1} = \lambda_1^{-1} + \lambda_3^{-1}.$
(In verifying this, the key relation is $p_2 = p_1 + p_3.)$
Hence we see that one necessary condition for the convexity
of~(\ref{x.92}) is
$\lambda_2^{-1} \leq \lambda_1^{-1} + \lambda_3^{-1};$ by symmetry,
the remaining conditions are $\lambda_1^{-1} \leq
\lambda_2^{-1} + \lambda_3^{-1}$ and
$\lambda_3^{-1} \leq \lambda_1^{-1} + \lambda_2^{-1}.$
It is easy to verify that if the $\lambda_i$ are also allowed to be
zero, and we write $0^{-1}=\infty$ and consider $\infty$ greater
than all real numbers, then these three conditions
continue to be necessary and sufficient for convexity.

Let us therefore index elements (\ref{x.92}) by the three parameters
$\lambda_1^{-1},$ $\lambda_2^{-1},$ $\lambda_3^{-1},$ defining
\[\label{x.94}
[a_1,\<a_2,\<a_3]\ =\ a_1^{-1}S_1\,\cup\,a_2^{-1}S_2\,\cup\,a_3^{-1}S_3
\quad (a_1, a_2, a_3\in[1,\infty]),\]
so that the lattice of centrally symmetric elements
of $\RelConv(S)_{\geq\{0\}}$ consists of the
sets $[a_1,a_2,a_3]$ with
\[\label{x.95}
a_1\,\leq\,a_2+a_3,\qquad a_2\,\leq\,a_1+a_3,\qquad
a_3\,\leq\,a_1+a_2.\]
Note that the ordering of this lattice is by reverse componentwise
comparison of expressions $[a_1,a_2,a_3],$ and
that $S_1,$ $S_2,$ $S_3$ are the
elements $[1,\infty,\infty],$ $[\infty,1,\infty],$ $[\infty,\infty,1].$
Lattice-theoretic meet is clearly given by componentwise
supremum, while the lattice-theoretic join of two elements is gotten
by first taking their componentwise infimum, which
represents their set-theoretic union, then getting its
relative convex hull by reducing the largest entry to the sum
of the other two if it exceeds this.
So, for instance, $S_1\vee S_2=
[1,\infty,\infty]\,\vee\,[\infty,1,\infty]$ is gotten by forming the
componentwise infimum, $[1,1,\infty],$ and then decreasing the last
component to the sum of the first two, getting $[1,1,2]$
(confirming the value $\lambda=1/2$ in our earlier description of this
element).
We can now calculate, e.g., the meet (componentwise supremum) of this
with $S_3=[\infty,\infty,1],$ namely $[\infty,\infty,2];$ and
take the join of this meet with $S_2=[\infty,1,\infty]$ by forming
the componentwise infimum $[\infty,1,2],$ and adjusting
the first component as above, getting $[3,1,2].$

Note that the only way these lattice operations yield components in
their values that did not occur as components in their
arguments is by addition; hence, as $L$ was defined to be generated
by $S_1,$ $S_2$ and $S_3,$ all
finite components $a_i$ that occur in the
expressions~(\ref{x.94}) for elements of $L$ are positive integers.
It is not hard to verify that all positive integers indeed occur,
and that the elements of $L$ are all the elements~(\ref{x.94})
satisfying~(\ref{x.95}) with
$a_1,\,a_2,\,a_3\in\{1, 2, 3, \dots, \infty\}.$
So we have a very arithmetic description of this lattice.

Though we have seen that $L$ contains an infinite descending chain,
it is interesting to note that the sublattice generated by any
finite set of elements~(\ref{x.94}),
none of which have any infinite components, is finite;
for the lattice operations will not produce, in any position, entries
larger than the corresponding entries of their arguments.\smallskip

It would be of interest to examine more general lattices of the
form $\RelConv(S)_{\geq\{0\}}$ for sets $S$ which
are unions of finitely many
line-segments (or rays) through the origin in $\R^n.$
In this situation, the conditions for convexity are always given by
linear inequalities in the ``$\lambda_i^{-1}\!$''; let me sketch why.

First, some general observations.
Suppose $p_1, \dots, p_m$ $(m\geq 3)$ are a minimal linearly dependent
family of vectors in $\R^n;$ thus they satisfy a nontrivial
linear relation $\sum c_i\,p_i=0,$ unique up to scalars.
For which families of positive real numbers
$\lambda_1,\dots,\lambda_m$ will the
points $\lambda_1\<p_1, \dots, \lambda_m\<p_m$ lie in
an \!$(m{-}1)$\!-dimensional affine subspace of their
\!$m$\!-dimensional span?
This will hold if and only if the linear relation
satisfied by these modified elements is an {\em affine relation}; i.e.,
has coefficients summing to $0.$
Those coefficients are $c_i\<\lambda_i^{-1},$ so the condition is
$$\textstyle\sum c_i\<\lambda_i^{-1}=0.$$
In this situation, can the affine relation among the
$\lambda_i\<p_i$ be written as an expression for one of them, say
$\lambda_j\<p_j,$ as a {\em convex} linear combination of the others?
One sees that this is so if and only if the coefficient
$c_j$ is opposite in sign to all the other $c_i;$ in that situation,
let us rewrite the expression satisfied by the $\lambda_i^{-1}$ as
$$\textstyle
\lambda_j^{-1}\ =\ -c_j^{-1}\sum_{i\neq j}\,c_i\<\lambda_i^{-1}.$$
Note that the abovementioned condition on signs can be satisfied by
at most one $j,$ and can be looked at as saying that the ray determined
by the corresponding vector $p_j$ is in the convex hull of the
rays determined by the other $p_i.$
When it is satisfied, one finds that a union of intervals
$\bigcup_i\,[0,\,\lambda_i\<p_i]$ is relatively convex in
the union of the rays determined by the $p_i$ if and only if
$\lambda_j$ has at least the value given by the above formula, i.e.,
if and only if
\[\label{x.99}
\lambda_j^{-1}\ \leq\ -c_j^{-1}\sum_{i\neq j}\,c_i\<\lambda_i^{-1}.\]
If the $c_i$ do not consist of one of one
sign and the rest of the opposite sign, then none of the
rays determined by the $p_i$ is in the convex hull of the rest, and
every set $\bigcup_i\,[0,\,\lambda_i\<p_i]$
$(\lambda_1,\dots,\lambda_m > 0)$ is relatively convex.
(An example of this situation is given for $n=3$ by letting
$p_1,\dots,p_4$ be the vertices of
a convex quadrilateral lying in a plane not containing $0.)$

Let us now drop the condition that the $p_i$ are minimal among linearly
dependent families, assuming only that none of them is a nonnegative
multiple of another (i.e., that they determine distinct rays), and
let $S$ denote the union of the rays through $0$ that they determine.
Then one can show that a union of intervals
$\bigcup_i\,[0,\,\lambda_i\<p_i]$ is relatively convex in
$S$ if and only if~(\ref{x.99}) holds for each
minimal linearly dependent family of $\geq 3$ of the $p_i$ whose
unique linear relation has exactly one coefficient $c_j$ of different
sign from the rest.
(The reduction to the minimal-linearly-independent-family case can
be gotten by a recursive application of Carath\'eodory's Theorem,
with $p_0\,{=}\,0,$ within subspaces spanned by successively
smaller subsets of $\{p_1,\,\dots\,,\,p_n\}.)$

As an example of the sort of lattice one gets,
let us drop the central-symmetry condition from
our snowflake construction, writing $-p_1,\,-p_2,\,-p_3$ as
$p_4,\,p_5,\,p_6,$ and considering general sets
$$\textstyle
\lambda_1\,\ch(0,p_1)\,\cup\,\dots\,\cup\,\lambda_6\,\ch(0,p_6)\quad
(\lambda_1,\dots,\lambda_6\in (0,\infty]).$$
Then setting $a_i=\lambda_i^{-1},$ we
get six conditions for relative convexity, namely
$$\textstyle a_2 \leq a_1+a_3,\quad a_3 \leq a_2+a_4,\quad
\dots,\quad a_1 \leq a_6+a_2,$$
each corresponding to the fact that one of our six rays lies
in the cone spanned by its two immediate neighbors.\smallskip

Our snowflake example showed that a
lattice $\RelConv(S)_{\geq\{0\}}$ could contain
a \!$3$\!-generator sublattice with an infinite descending chain.
Let me sketch an example with an infinite ascending chain.
In $\R^2,$ let
$$\textstyle p_1=(0,3),\quad p_2=(1,2),\quad
p_3=(2,1),\quad p_4=(3,0).$$
(Or for a more abstract description, take any four points
of $\R^n$ in arithmetic progression, on a line
not passing through the origin.)
Define
$$\textstyle S_i=\ch(0,p_i)\quad\text{and}
\quad S=S_1\cup S_2\cup S_3\cup S_4.$$
Now let $L$ be the sublattice of $\RelConv(S)_{\geq\{0\}}$ generated by
$$\textstyle x_1=S_1\cup(S_2/2),\quad y=S_2\cup S_3,\quad
x_2=(S_3/2)\cup S_4.$$
It is easy to see from a sketch that,
starting with $x_1\wedge y,$ if we alternately
apply $(-\vee x_2)\wedge y$ and $(-\vee x_1)\wedge y,$ we
obtain an infinite ascending chain of subsets.

Is the existence of infinite chains in sublattices
of lattices $\RelConv(S)$ generated by few elements limited to cases
where $S\neq \R^n,$ or does it also occur in
lattices $\Conv(\R^n)?$
To get a large part of the answer without any computation, recall
that $\RelConv(S)$ is isomorphic to the lattice
of convex hulls in $\R^n$ of subsets of $S.$
This lattice, which we shall here
denote $L_S,$ is a subset but not a
sublattice of $\Conv(\R^n);$ however, in cases
like those considered above, where $S$ is a finite union of
convex sets, $S=S_1\cup \dots\cup S_m,$ we
can write the operations of this lattice as
``polynomial operations'' in those of $\Conv(\R^n).$
Namely, temporarily writing $\wedge_S$ and $\vee_S$ for
the operations of $L_S,$ and $\wedge$ and $\vee$ for
those of $\Conv(\R^n),$ we get, for $x, y\in L_S,$
\[\label{x.96}
x\vee_S y=x\vee y,\qquad x\wedge_S
y=\bigvee_{i}\,(x\wedge y\wedge S_i).\]
Hence if the sublattice of $L_S$ generated
by elements $y_1,\dots,y_k$ has
an infinite ascending or descending chain (or any other
specified join-sublattice), so will the sublattice
of $\Conv(\R^n)$ generated by $y_1, \dots, y_k, S_1, \dots, S_m.$
In the case of our ``snowflake lattice'', the $y$'s and
the $S$'s happen to be the same, so we immediately conclude that
the sublattice of $\Conv(\R^2)_{\geq\{0\}}$ generated
by those three elements has an infinite descending chain.
From our example with an ascending chain, the best conclusion this
general argument gives is that the sublattice
of $\Conv(\R^2)_{\geq\{0\}}$ generated by the six
elements $S_1,$ $S_2,$ $S_3,$ $S_4,$ $S_2/2$ and $S_3/2$ has an
infinite ascending chain.
Nevertheless, a little diagram-drawing shows that in this case,
the sublattice of $\Conv(\R^2)_{\geq\{0\}}$ generated by the
three elements $\ch(x_1),$ $\ch(x_2)$ and $\ch(y)$ shows essentially
the same behavior as our lattice of relatively convex subsets.

I do not know any examples of \!$3$\!-generator lattices of convex sets
(relative or absolute) that have both infinite ascending and
descending chains, or that have infinite antichains.
On the other hand, it is not hard to show that the \!$4$\!-generator
sublattice of $\Conv(\R^2)$ generated by the diameters of a
regular octagon has all three.

\section{Notes on related work on relatively convex sets,
and some further observations.}\label{NHrlmsc}
The lattices $\RelConv(S)$ are examples of what are known as
convex geometries; for the definition,
and results on these, see \cite{PE+RJ}, \cite{KA.CvG}.\smallskip

Huhn \cite{AH.n-d} looked briefly at lattices of relatively
convex sets determined by finite sets $S \subseteq\R^n$ for
the purpose of ``approximating'' $\Conv(\R^n)$ by
finite lattices, and for the same purpose he considered
in \cite[\S2]{AH.Dop} the dual construction, namely the lattice of
those subsets of $\R^n$ which can be represented as intersections of
members of a given finite set of closed half-spaces.
Not surprisingly in view of the dual natures of these two sorts of
relativization, he found that lattices of the latter sort satisfied the
identities of the form $D_n\op$ that he had obtained in
the nonrelativized lattice, but not those of the
form $D_n$ (cf.\ Proposition~\ref{X.88} and Lemma~\ref{X.90} above).

In \cite{FW+MS} it is shown that for every finite
lattice $L$ which is ``lower bounded'' (a strengthening of join
semidistributive), there exist an $n$ and a finite
subset $S \subseteq \mathbb Q^n$ such that $L$ is
embeddable in the sublattice of $\Conv(\R^n)$ generated
by $\{\{p\}\mid p\in S\}.$
Clearly, such an embedding will send all members of $L$ to convex
polytopes; if we let $S'$ denote the union of the vertex-sets of this
finite set of polytopes, it is not hard to see that we get an embedding
of $L$ in $\RelConv(S').$
Whether embeddability in the lattice of relatively convex sets of
a finite subset of $\R^n$ holds not only for lower bounded lattices, but
for all join-semidistributive lattices, is an open
question \cite[Problem~1]{FW+MS}, cf.\ \cite[Problem~3]{KA.CvG}.
It is also shown in \cite{FW+MS} that every lattice can be
embedded in the lattice of convex subsets of some infinite-dimensional
vector space (over an arbitrary totally ordered division ring).
In particular, such lattices need not satisfy any
nontrivial lattice identities.
In contrast, the lattice of {\em subspaces} of any vector space
satisfies the modular identity and others.

Here is another embedding result using not necessarily
finite-dimensional vector spaces, although it may be seen that the
connection with convexity is somewhat artificial, based on the fact
that subspaces are in particular convex sets; it is essentially
a result on lattices of ``relative subspaces''.
(Since most of the results in this note concerned subspaces
of $\R^n,$ we only gave our general definitions for that
case, but we shall use them here without that restriction.)
\begin{lemma}\label{X.98}
Let $V$ be a real vector space and $B$ a basis of $V.$
For each pair of distinct elements $a,b\in B,$ let $S_{a,b}$ be
the subspace of $V$ spanned by $a-b.$
Define the set
$$\textstyle S\ =\ \bigcup_{a,b}\,S_{a,b}\ \subseteq\ V.$$
Then the following lattices are isomorphic:

{\rm(a)} The lattice $\mathrm{Equiv}(B)$ of all
equivalence relations on $B$ \r(ordered by inclusion\r).

{\rm(b)} The lattice $\mathrm{RelSubSp}(S)$ consisting
of all sets of the form $S\cap U$ where $U$ is a subspace of $V.$
\r(Here meet is given by intersection; join by
taking the subspace spanned by the union of the sets in
question, and intersecting this with $S.)$

{\rm(b$\!{}'$)} The lattice of all subspaces of $V$ spanned
by subsets of $S.$
\r(The join of two such subspaces $W_1,$ $W_2$ is $W_1+W_2,$ the
meet is the span of $W_1\ \cap\ W_2\ \cap\ S.)$

{\rm(c)} The sublattice of $\RelConv(S)$ consisting
of all elements thereof that are unions of subspaces~$S_{a,b}.$

{\rm(c$\!{}'$)} The lattice of all subsets of $V$ which are convex
hulls of unions of subpaces $S_{a,b}.$
\r(Join as in $\Conv(\R^n),$ while
the meet of $x_1,$ $x_2$ is $\ch(x_1\cap x_2\cap S)).$
\end{lemma}
\begin{proof}
That the sets described in (b), (b$\!{}')$ and (c$\!{}'),$ ordered by
inclusion, form lattices, with meet and join as described, is immediate.

Note that the convex hull of a union of subspaces of $V$ is the
sum of those subspaces, and that the intersection
of $S$ with a subspace is always a union of certain
of the $S_{a,b}.$
From this it is easily seen that~(b) and~(c) are not merely, as
asserted, isomorphic, but equal, and likewise~(b$\!{}')$ and~(c$\!{}'$).
It is also clear that~(b) is isomorphic to~(b$\!{}'$), via the ``span
of'' map in one direction and the operator $-\cap S$ in
the other.
So these four lattices are isomorphic; to complete the
proof we shall describe an isomorphism between the
lattice $\mathrm{Equiv}(B)$ of~(a) and the lattice of~(b$\!{}'$).

Given an equivalence
relation $R\in \mathrm{Equiv}(B),$ let $\phi(R)$ be
the subspace of $V$ spanned by all
elements $a-b$ with $(a,b)\in R,$ which
by definition belongs to~(b$\!{}'$), while
given a subspace $W \subseteq V$ spanned by a subset
of $S,$ let $\psi(W) = \{(a,b)\mid a-b\in W\},$ which it is
easy to check is an equivalence relation on $B.$
The maps $\phi$ and $\psi$ are clearly isotone, and from the
definition of~(b$\!{}'$), we see
that $\phi\psi$ is the identity function thereof;
moreover, for any $R\in \mathrm{Equiv}(B)$ it is clear
that $\psi(\phi(R)) \geq R,$ so it remains
to prove the reverse inequality.

So suppose that $(a,b)\in\psi(\phi(R)),$ i.e.,
that $a-b\in\phi(R).$
From the definition of $\phi(R)$ it is easy to see that
for each \!$R$\!-equivalence class $C \subseteq B,$ the
sum of the coefficients of all members of $C$ in any element
of $\phi(R)$ is zero.
But the only way this can hold for $a-b$ is
if $a$ and $b$ are in the same equivalence class,
i.e., $(a,b)\in R,$ as required to complete
our proof.\end{proof}

Pudl\'ak and T\accent 23uma \cite{PT.lat}
have shown that every finite lattice $L$ embeds
in $\mathrm{Equiv}(X)$ for some finite set $X.$
Hence by the above lemma, for every such $L$ one
can find an $n$ and a subset $S \subseteq \R^n$ such
that $L$ is embeddable in $\RelConv(S).$
It would be interesting to know whether this same conclusion can be
proved without using the deep result of \cite{PT.lat}.
So far as is known, the embedding of $L$ in a
lattice $\mathrm{Equiv}(X)$ may require a
set $X$ whose cardinality is enormous compared with that
of $L$ (see \cite{HK.lat}, \cite{TI.cong} for some improvements
on the bound of \cite{PT.lat}), but it is plausible that one could do
better with embeddings in lattices of relatively convex sets.\smallskip

Our final remark will concern the observations at the
end of the preceding section on the form of the lattice operations
of $\RelConv(S)$ when $S$ is a finite union of
convex sets $S_1\cup \dots\cup S_m.$
Let us put these in a more general context.
(Readers allergic to category theory may ignore this discussion.)

Let $\Lat$ denote the category of all
lattices, objects of which we will here
write $L=(|L|,\ \vee,\ \wedge),$ distinguishing
between the lattice $L$ and its underlying set $|L|.$
For $m$ a natural number, let $\Latm$ denote the category of lattices
with $m$ distinguished elements,
i.e., of systems $(|L|, \vee, \wedge, S_1,\dots,S_m)$ such
that $(|L|,\ \vee,\ \wedge)$ is a lattice and $S_1,\dots,S_m\in |L|,$
and where a morphism between such systems means
a lattice homomorphism which respects the ordered \!$m$\!-tuple
of distinguished elements.
Let us, finally, write $\Lat^{\inr}$ for the category
of objects $(|L|,\ \vee,\ \wedge,\ \inr),$ where
$(|L|,\ \vee,\ \wedge)$ is a lattice,
and $\inr$ is an ``interior operator'' (the dual of a closure operator),
that is, a map $|L| \to |L|$ satisfying
$$\textstyle\inr(x)\leq x,\quad x\leq y\implies \inr(x)\leq\inr(y),\quad
\inr(\inr(x))=\inr(x).$$
(Here by ``$u\leq v\!$'' we of course
mean $u=u\wedge v,$ equivalently, $u\vee v=v.)$
The morphisms of $\Lat^{\inr}$ will be the
lattice homomorphisms respecting this additional operation.

We can define a functor $\Latm\to \Lat^{\inr}$ taking
each object $(|L|,\vee, \wedge,\linebreak[2] S_1,\dots,S_m)$ to
the object $(|L|,\vee,\wedge,\inr_{(S_i)}),$ where
the interior operator is defined by
$$\textstyle \inr_{(S_i)}(x)\ =\ \bigvee_i (x\wedge S_i),$$
and another functor $\Lat^{\inr} \to \Lat,$ taking
each object $(|L|,\ \vee,\ \wedge,\ \inr)$ to the object
$(|L|^{\inr},\ \vee,\ \wedge^{\inr}),$ where
$$\textstyle |L|^{\inr}=\{x\in|L|\mid x=\inr(x)\}
\text{\quad and\quad}x\,\wedge^{\inr}\,y\ =\ \inr(x\wedge y).$$
(It is not hard to verify that the join operation of $L$
carries $|L|^\inr$ into itself.)

We now see that if we take $L=\Conv(\R^n),$ and
let $S_1,\dots,S_m$ be any $m$ elements of this lattice,
then the composite of the above two functors, applied
to $(|L|,\ \vee,\ \wedge,\ S_1,\dots,S_m),$ gives precisely the
lattice we named $L_S\cong\RelConv(S),$ for $S=S_1\cup\dots\cup S_m.$

The constructions given by the above functor
$\Lat^{\inr} \to\Lat,$ and its dual, with a closure operator replacing
the interior operator, are well-known, if not in this functorial form.
I do not know whether the
construction $\Latm \to\Lat^{\inr}$ (and its variant
with lattices replaced by complete lattices and the
specified finite family by an arbitrary family) has been considered.

Stepping back a little further, we may observe that the
lattice $\Conv(\R^n)$ arises as the fixed set of the
closure operator $\ch(-)$ on the lattice $2^{\R^n}$ of
subsets of $\R^n,$ so that the construction of the
mutually isomorphic lattices $L_S$ and $\RelConv(S)$ can
be seen as arising from the interaction of the closure
operator $\ch(-)$ and the interior
operator $-\cap S$ on $2^{\R^n}.$
Again, this situation can be made into a general construction.

The reader familiar with the concept of representable
algebra-valued functors (\cite[Chapter~9]{GB.245}
or \cite[\S1, \S8]{coalg}) will be happy to observe that all
the functors of the above discussion are representable.

\section{A question.}\label{NHRefQ}
The referee has pointed out that some properties of the lattice of
convex sets are known to change if the base field $\R$ is replaced by
another ordered field (e.g., $\mathbb Q\<),$ but that the
arguments of \S\S\ref{NHintro}-\ref{NH_Radon}
look as though they should work over any
ordered field; perhaps even any ordered division ring.
I have the same feeling, but as as an amateur in the area,
I will leave this question to others.
(I do not know whether the theorems of Helly, Carath\'eodory and Radon
hold in that context, nor how much of what I have justified
as geometrically evident may rely on properties of the reals.)

A straightforward generalization of these results could not, of course,
extend to \S\ref{NHcp}, on compact convex sets, since over an ordered
field which is not locally compact, the only nonempty compact convex
sets are the singletons, which do not form a lattice.
One might be able to prove results like those of that section
with ``compact'' replaced by ``closed and bounded'', but new proofs
would be needed, since the theorem of Minkowski we used
there is not true in that context.
Later sections depend to varying degrees on that one.

\end{document}